\documentclass[review]{elsarticle}          

\usepackage[english]{babel}
\usepackage[utf8]{inputenc}
\usepackage[T1]{fontenc}
\usepackage{url}

\usepackage{amssymb,amsmath,amscd}

\catcode`@=11
\def\plslash{\ifx\@currsize\normalsize
{\mathchoice
{\mbox{\raisebox{0.2ex}{$\scriptstyle\circ$}\kern-1ex$\setminus$}}
{\mbox{\raisebox{0.2ex}{$\scriptstyle\circ$}\kern-1ex$\setminus$}}%
{\mbox{\raisebox{0.14ex}{$\scriptscriptstyle\circ$}\kern-0.8ex%
${\scriptstyle\setminus}$}}%
{\,\mbox{\raisebox{0.14ex}{$\scriptscriptstyle\circ$}\kern-0.8ex%
${\scriptstyle\setminus}$}}}%
\else\ifx\@currsize\large\,\mbox{\raisebox{0.2ex}{$\scriptstyle\circ$}\kern-1ex$\setminus$}
\else\ifx\@currsize\small\,\mbox{\raisebox{0.2ex}{$\scriptstyle\circ$}\kern-1ex$\setminus$}
\else\,\mbox{\raisebox{0.2ex}{$\scriptstyle\circ$}\kern-1ex$\setminus$}
\fi\fi\fi}

\def\prslash{\ifx\@currsize\normalsize
{\mathchoice
{\mbox{\raisebox{0.2ex}{$\scriptstyle\circ$}\kern-1ex$/$}}
{\mbox{\raisebox{0.2ex}{$\scriptstyle\circ$}\kern-1ex$/$}}%
{\mbox{\raisebox{0.14ex}{$\scriptscriptstyle\circ$}\kern-0.8ex%
${\scriptstyle/}$}}%
{\mbox{\raisebox{0.14ex}{$\scriptscriptstyle\circ$}\kern-0.8ex%
${\scriptstyle/}$}}}%
\else\ifx\@currsize\large\mbox{\raisebox{0.2ex}{$\scriptstyle\circ$}\kern-1ex$/$}
\else\ifx\@currsize\small\mbox{\raisebox{0.2ex}{$\scriptstyle\circ$}\kern-1ex$/$}
\else\mbox{\raisebox{0.2ex}{$\scriptstyle\circ$}\kern-1ex$/$}
\fi\fi\fi}

\newcommand{\lslash}{\protect\plslash}
\newcommand{\rslash}{\protect\prslash}
\def\pomoins{\ifx\@currsize\normalsize
\mbox{ $\circ\kern-1.48ex-$ }
\else\ifx\@currsize\large\mbox{ $\circ\kern-1.45ex-$ }
\else\ifx\@currsize\small\mbox{ $\circ\kern-1.51ex-$ }
\else\mbox{ $\circ\kern-1.40ex-$ }\fi\fi\fi}

\catcode`@=12

\newcommand{\rfrac}[2]{\,\rule[0.2ex]{0.1ex}{0.4ex}\kern-0.27ex%
\frac{#1\,}{\,#2}\kern-0.29ex\rule[0.55ex]{0.1ex}{0.4ex}\,\mbox{}}

\newcommand{\lfrac}[2]{\,\rule[0.55ex]{0.1ex}{0.4ex}\kern-0.29ex%
\frac{\,#1}{#2\,}\kern-0.27ex\rule[0.2ex]{0.1ex}{0.4ex}\,\mbox{}}

\catcode`@=11
\def\plslashh{\ifx\@currsize\normalsize
{\mathchoice
{\mbox{\raisebox{0.2ex}{$\scriptstyle\bullet$}\kern-1ex$\setminus$}}
{\mbox{\raisebox{0.2ex}{$\scriptstyle\bullet$}\kern-1ex$\setminus$}}%
{\mbox{\raisebox{0.14ex}{$\scriptscriptstyle\bullet$}\kern-0.8ex%
${\scriptstyle\setminus}$}}%
{\,\mbox{\raisebox{0.14ex}{$\scriptscriptstyle\bullet$}\kern-0.8ex%
${\scriptstyle\setminus}$}}}%
\else\ifx\@currsize\large\,\mbox{\raisebox{0.2ex}{$\scriptstyle\bullet$}\kern-1ex$\setminus$}
\else\ifx\@currsize\small\,\mbox{\raisebox{0.2ex}{$\scriptstyle\bullet$}\kern-1ex$\setminus$}
\else\,\mbox{\raisebox{0.2ex}{$\scriptstyle\bullet$}\kern-1ex$\setminus$}
\fi\fi\fi}

\def\prslashh{\ifx\@currsize\normalsize
{\mathchoice
{\mbox{\raisebox{0.2ex}{$\scriptstyle\bullet$}\kern-1ex$/$}}
{\mbox{\raisebox{0.2ex}{$\scriptstyle\bullet$}\kern-1ex$/$}}%
{\mbox{\raisebox{0.14ex}{$\scriptscriptstyle\bullet$}\kern-0.8ex%
${\scriptstyle/}$}}%
{\mbox{\raisebox{0.14ex}{$\scriptscriptstyle\bullet$}\kern-0.8ex%
${\scriptstyle/}$}}}%
\else\ifx\@currsize\large\mbox{\raisebox{0.2ex}{$\scriptstyle\bullet$}\kern-1ex$/$}
\else\ifx\@currsize\small\mbox{\raisebox{0.2ex}{$\scriptstyle\bullet$}\kern-1ex$/$}
\else\mbox{\raisebox{0.2ex}{$\scriptstyle\bullet$}\kern-1ex$/$}
\fi\fi\fi}

\newcommand{\lslashh}{\protect\plslashh}

\newcommand{\id}{\mathsf{Id}}
\newcommand{\im}{\mathsf{Im}}
\newcommand{\zmax}{\overline{\mathbb{Z}}_{\max}}

\newcommand{\zmaxgamma}{\gamma^\ast\overline{\mathbb{Z}}_{\max}[\![\gamma]\!]}

\newcommand{\cod}{\mathcal{C}_{\mathtt{O}}(\mathcal{S})}
\newcommand{\cd}{\mathcal{C}(\mathcal{S})}
\newcommand{\Zset}{\mathbb{Z}}
\newcommand{\Nset}{\mathbb{N}_0}
\newcommand{\sd}{\mathcal{S}}
\newcommand{\Is}{\mathrm{I}\mathcal{S}}
\newcommand{\sseries}{\mathcal{S}[\![ z_1,...,z_p]\!]}
\newtheorem{thm}{Theorem}
\newtheorem{propriete}[thm]{Property}
\newtheorem{proposition}[thm]{Proposition}
\newtheorem{cor}[thm]{Corollary}
\newtheorem{notation}[thm]{Notation}
\newtheorem{defn}[thm]{Definition}
\newtheorem{exmp}[thm]{Example}
\newproof{pf}{Proof}
\newtheorem{lem}[thm]{Lemma}
\newtheorem{rem}[thm]{Remark}

\begin{document}

\begin{frontmatter}

\title{Duality and interval analysis over idempotent semirings\footnote{This work was supported by the french/german project Procope/DAAD  21994UH, and also by the french/brazilian project CNPq, FAPEMIG, CAPES/COFECUB  Ph 642/09.}
\small{Extended version of 
Linear Algebra and its Application, 437 (2012) 2436–2454  }
} 
 \author[label1,label2]{Thomas Brunsch}
 \author[label1]{Laurent Hardouin\corref{cor1}}
 \author[label3]{Carlos Andrey Maia}   
 \author[label2,label4]{J\"org Raisch}

\address[label1]{Laboratoire d'Ing\'enierie des Syst\`emes Automatis\'es, \\Universit\'e d'Angers,\\ 62, avenue Notre Dame du Lac, 49000 Angers, France.}  
\address[label2]{Fachgebiet Regelungssysteme, Technische Universit\"at Berlin, Einsteinufer 17, 10587 Berlin, Germany}
\address[label3]{ Departamento de Engenharia El\'{e}trica, Universidade Federal de Minas
Gerais (UFMG). Av. Ant\^{o}nio Carlos 6627, Pampulha, 31270-010, Belo
Horizonte, MG, Brazil.}
\address[label4]{Fachgruppe System- und Regelungstheorie, Max-Planck-Institut f\"ur Dynamik komplexer technischer Systeme, Sandtorstr. 1, 39106 Magdeburg, Germany}
\cortext[cor1]{Corresponding author : \url{laurent.hardouin@univ-angers.fr}}


\begin{abstract}
In this paper semirings with an idempotent addition are considered. These
algebraic structures are endowed with a partial order. This allows to consider residuated maps to
solve systems of inequalities $A \otimes X \preceq B$ (see \cite{blyth72}).
The purpose of this paper is to consider a dual product, denoted $\odot$, and the dual residuation of matrices,
in order to solve the following inequality $ A \otimes  X \preceq X  \preceq B \odot X$.
Sufficient conditions ensuring the existence of a non-linear projector in the solution set are proposed.
The results are extended to semirings of intervals such as they were introduced in \cite{litvinov01a}.

\end{abstract}

\end{frontmatter}

\section{Introduction}
Many problems in mathematics are non-linear in the traditional sense but appear to
be linear over idempotent semirings. The max-plus algebra is a popular semiring widely studied (see e.g., \cite{carre79,cuning,butkovic07,butkovic09,Lorenzoa11}).
An idempotent semiring $\mathcal{S}$ can be endowed with a partial order relation. According to this order relation, and according to continuity assumptions, it is possible to obtain the greatest solution of inequality $A \otimes X\preceq B$ where $A,X$ and $B$ are matrices of proper dimension and $(A\otimes X)_{ij}=\bigoplus \limits_{k=1\ldots n} \left ( a_{ik}\otimes x_{kj} \right )$. The greatest solution is obtained by considering residuation theory.
In this paper we will consider the dual matrix product $A \odot X$ defined as $(A \odot X)_{ij}=\bigwedge \limits_{k=1\ldots n} \left ( a_{ik}\odot x_{kj} \right )$, where $\wedge$ represents the greatest lower bound. Then we will consider the dual residuation to deal with computation of the smallest solution of inequality $ A \odot X \succeq B$. The existence of a unique solution is not always ensured.
 Nevertheless if all elements of the semiring admit an inverse ($i.e.$,  it is a semifield) then the smallest solution exists. This condition is fulfilled in (max-plus) algebra and it has allowed to deal with opposite semimodules in \cite{cohen04}. This condition is  fulfilled neither in the semirings of non decreasing power series nor in  the semirings of intervals such as introduced in \cite{litvinov01c,lhommeau04,hardouin09b}, hence we will give  some sufficient conditions to ensure the existence of this smallest solution.

From a practical point of view, it is useful to be able to solve systems such as $A \otimes X \preceq X \preceq B \odot X$, as they are involved in the study of dynamical discrete event systems subject to constraints (see \cite{ouerghi06b, brunsch11,houssin06b}).
Hence sufficient conditions for the existence of a projector in the set of solutions is given. Its
computation is based on additive closure of matrices and on the dual residuation of the dual product.
This projector is also given in semirings of intervals which allow us to deal with uncertainties.

This paper is organized as follows:
in Section \ref{preliminaries}, algebraic preliminaries are recalled. More precisely, semiring defintion are first introduced and then some useful theorems about residuation theory are recalled. Next the section is devoted to the presentation of closure mapping properties.
In Section \ref{section:dualproduct}, the dual product  and its dual residuation are considered.
Inequalities $ A\otimes X \preceq X \preceq B \odot X$ is considered in Section \ref{section:axbx}, and in order to propose a projector in the solution set, some sufficient conditions are given. In Section \ref{section:Example}, the previous results are applied in the semiring (max,plus) and in a semiring of non-decreasing power series.
In Section \ref{section:intervalanalysis}, semirings of intervals are considered. Useful results initially presented in \cite{litvinov01c,lhommeau04,hardouin09b} are recalled, and the results of Section \ref{section:dualproduct} are extended in this algebraic setting.

\section{Preliminaries} \label{preliminaries}
\subsection{Idempotent Semiring}
In this section we recall useful results (for a more exhaustive presentation see reference   \cite{BICOQ}).
\begin{defn}[Monoid]
$(M,\cdot,e)$ is a monoid if $\cdot$ is an internal law, associative and with an identity element $e$. If the law
$\cdot$ is commutative, $(M,\cdot,e)$ is a commutative monoid.
\end{defn}
\begin{defn}[Idempotent Semiring, semifield] \label{def:semiring} An idempotent \textit{semiring} is a set,  $\mathcal{S}$,  endowed with two
internal operations denoted by $\oplus$ (addition) and $\otimes$ (multiplication) such that :
\begin{itemize}
\item[] ($\sd,\oplus,\varepsilon$) is an idempotent commutative monoid, $i.e.$, $\forall a \in \sd, a\oplus a =a$,
\item[] ($\sd,\otimes,e$) is a monoid,
\item[] $\otimes$ operation is distributive with respect to $\oplus$,
\item[] $\varepsilon$ is absorbing for the law $\otimes$, $i.e., \text{ } \forall a,~~\varepsilon \otimes a = a
\otimes \varepsilon = \varepsilon$.
\item[] If $\otimes$ is commutative, the semiring is said to be commutative.
A \textit{semifield} is a semiring in which all elements except $\varepsilon$ have a multiplicative inverse.
\end{itemize}
\end{defn}

An idempotent semiring\footnote{In the following we will only refer to idempotent semirings and therefore drop the adjective} can be endowed with a canonical \textit{order} defined  by: $ a
\succeq b~$ iff $~a=a \oplus b$. Then it becomes a sup-semilattice,
and $a \oplus b$ is the least upper bound of $a$ and $b$. A
semiring is \textit{complete} if sums of infinite number of terms
are always defined, and if multiplication distributes over
infinite sums, too. In particular, the sum of all elements of a complete
semiring is defined and denoted by $\top$ (for "top"). A complete
semiring  becomes a complete lattice for which
the greatest lower bound of $a$ and $b$ is denoted $a \wedge b$.


\begin{defn}[Subsemiring]
\label{def:subsemiring} A subset $\mathcal{C}$ of a semiring is
called a subsemiring of $\mathcal{S}$ if
\begin{enumerate}
\item[] $\varepsilon \in \mathcal{C}$ and $e \in
\mathcal{C}$ ;
\item[] $\mathcal{C}$ is closed for
$\oplus$ and $\otimes$, i.e, $\forall a,b \in \mathcal{C}$, $a
\oplus b \in \mathcal{C}$ and $a \otimes b \in \mathcal{C}$.
\end{enumerate}
Furthermore the subsemiring is complete if it is closed for infinite sums and if the product distributes over infinite sums.
\end{defn}

\begin{lem}[{\citep[$\S 4.3.4$]{BICOQ}}]\label{lem:productdistributeoverinf}
Let $\sd$ be a semiring. $\forall a,b,c \in \mathcal{S}$ the following inequality holds :
$$\nonumber c \otimes (a\wedge b) \preceq (c \otimes a) \wedge (c \otimes b).$$
 Furthermore, if $c$ admits a multiplicative inverse, $i.e.$, if there exists a unique element, denoted $c^{-1}$, such that $c^{-1}\otimes c=c\otimes c^{-1}=e$, then $$ c \otimes (a\wedge b) = (c \otimes a) \wedge (c \otimes b) .$$
\end{lem}
%
%

\begin{defn}[Formal power series] A formal power series in $p$
(commutative) variables, denoted $z_1$ to $z_p$, with coefficients
in a semiring $\sd$, is a mapping $s$ defined from $\Zset^p$ into  $\sd$: $\forall k=(k_1,...,k_p) \in \Zset^p$,
$s(k)$ represents the coefficient of $z_1^{k_1}...z_p^{k_p}$ and
$(k_1,...,k_p)$ are the exponents. Another equivalent
representation is
$$ s(z_1,...,z_p )=\bigoplus_{k \in \Zset^p} s(k)z_1^{k_1}...z_p^{k_p}.
$$
\end{defn}

\begin{defn}[Semiring of series] \label{def:formalpowerseries}
The set of formal power series with coefficients in a semiring $\sd$ endowed with the following sum and Cauchy
product:
\begin{eqnarray}\label{operationseries}
\nonumber  s \oplus s'  : (s \oplus s')(k)= s(k) \oplus s'(k),\\
\nonumber  s \otimes s' : (s \otimes s') (k)=\bigoplus_{i+j=k} s(i)\otimes   s'(j),
\end{eqnarray}
is a semiring denoted $\sseries$. If $\mathcal{S}$ is complete, $\sseries$ is complete.
A series with a finite support is called a polynomial, and a monomial if there is only one element in the series.
The greatest lower bound of series is given by :
\begin{equation}
\nonumber s \wedge s'  : (s \wedge s')(k)= s(k) \wedge s'(k).
\end{equation}
\end{defn}


\subsection{Residuation Theory}
Residuation theory allows to deal with the inverse of order preserving mappings defined over
ordered sets, $i.e.$ a set equipped with a partial
order relation. This theory gives another point of view on Galois connection. Useful references are
\cite{croisot56,blyth72,blyth05}.

\begin{defn}[Continuity] \label{def_continuity}
An order preserving mapping $f : \mathcal{D} \rightarrow
\mathcal{E}$, where $\mathcal{D}$ and $\mathcal{E}$ are complete ordered
sets, is a mapping such that: $x\succeq y \Rightarrow f(x) \succeq f(y)$.
It is said to be isotone in \cite{BICOQ}.

A mapping $f$ is lower-semicontinuous ($l.s.c.$), respectively, upper-semicontinuous ($u.s.c.$)
if, for every (finite or infinite) subset $\mathcal{X}$ of $\mathcal{D}$,
$$ f(\bigoplus_{x \in \mathcal{X}} x)= \bigoplus_{x \in \mathcal{X}} f(x),$$
respectively,
$$ f(\bigwedge_{x \in \mathcal{X}} x)= \bigwedge_{x \in \mathcal{X}} f(x).$$
A mapping  $f$ is continuous if it is both $l.s.c.$ and $u.c.s$.
\end{defn}

\begin{defn}[Image, Kernel] \label{def:imagekernel}
Let $f : \mathcal{D} \rightarrow
\mathcal{E}$ be a mapping, where $\mathcal{D}$ and $\mathcal{E}$ are semirings.
The image of $f$, denoted $\im f$,  is classically defined as $\im f=\{y \in \mathcal{E} | y=f(x) \text{ for some } x  \in \mathcal{D} \} $.
The equivalence kernel is defined as $\mathrm{ker} f := \{(x,x')\in \mathcal{D} \times \mathcal{D} \mid f(x) = f(x')\}$.
\end{defn}

\begin{defn}[Residuated and dually residuated mapping] \label{def_residuation}

An order preserving mapping $f : \mathcal{D} \rightarrow
\mathcal{E}$, where $\mathcal{D}$ and $\mathcal{E}$ are ordered
sets, is a \textit{residuated mapping} if for all $y \in
\mathcal{E}$, the least upper bound of the subset $\{x | f(x)
\preceq y\}$ exists and belongs to this subset. It is then denoted
by $f^\sharp(y)$. The mapping $f^\sharp$ is called the residual of
$f$. When $f$ is residuated, $f^\sharp$ is the unique order
preserving mapping such that \small
\begin{equation}\label{residuated1}
\begin{array}{lcl}
f \circ f^\sharp \preceq \id_{\mathcal{E}} & \textnormal{~~~and~~~} &
f^\sharp \circ f \succeq \id_{\mathcal{D}},
\end{array}
\end{equation}
\normalsize where $\id$ is the identity mapping  on
$\mathcal{D}$ and $\mathcal{E}$ respectively.\\
Mapping $g$ is a \textit{dually residuated mapping} if for all $y \in
\mathcal{E}$, the greatest lower bound of the subset $\{x | g(x)
\succeq y\}$ exists and belongs to this subset. It is then denoted
by $g^\flat (y)$. The mapping $g^\flat$ is called the dual residual of
$g$. When $g$ is dually residuated, $g^\flat$ is the unique order
preserving mapping such that \small
\begin{equation}\label{dualresiduated1}
\begin{array}{lcl}
g \circ g^\flat \succeq \id_{\mathcal{E}} & \text{~~~and~~~} &
g^\flat \circ g \preceq \id_{\mathcal{D}}.
\end{array}
\end{equation}
\end{defn}

\begin{rem}\label{rem:gflatsharp=g}
According to this definition, it is clear that $f^\sharp$ is dually residuated and that $g^\flat$ is residuated,
furthermore, $(f^\sharp)^\flat=f$ and ($g^\flat)^\sharp=g$.
\end{rem}

\begin{thm}[{\citep[$\S 4.4.2$]{BICOQ}}]
\label{thm:residuation} Consider the order preserving mappings $f:\mathcal{E}
\rightarrow \mathcal{F}$ and $g:\mathcal{E}
\rightarrow \mathcal{F}$ where $\mathcal{E}$ and $\mathcal{F}$ are
complete semirings. Their bottom elements are, respectively,
denoted by $\varepsilon_{\mathcal{E}}$ and
$\varepsilon_{\mathcal{F}}$. Their top elements are, respectively,
denoted by $\top_{\mathcal{E}}$ and
$\top_{\mathcal{F}}$.\\ Mapping $f$ is residuated iff
$f(\varepsilon_{\mathcal{E}}) = \varepsilon_{\mathcal{F}}$ and
$f(\bigoplus_{x \in \mathcal{X}}x) = \bigoplus_{x \in
\mathcal{X}}f(x)$ for each $\mathcal{X} \subseteq \mathcal{E}$
(i.e., $f$ is \textit{lower-semicontinuous}), furthermore $f^\sharp(\top_{\mathcal{F}})=\top_{\mathcal{E}}$ and
$ f^\sharp(\bigwedge_{y \in \mathcal{Y}}y) = \bigwedge_{y \in
\mathcal{Y}}f^\sharp(y)$ for each $\mathcal{Y} \subseteq \mathcal{F}$ (i.e., $f^\sharp$ is  \textit{upper-semicontinuous}).\\
Mapping $g$ is dually residuated iff
$g(\top_{\mathcal{E}}) = \top_{\mathcal{F}}$ and
$g(\bigwedge_{x \in \mathcal{X}}x) = \bigwedge_{x \in
\mathcal{X}}g(x)$ for each $\mathcal{X} \subseteq \mathcal{E}$
(i.e., $g$ is \textit{upper-semicontinuous}), furthermore $g^\flat(\varepsilon_{\mathcal{F}}) = \varepsilon_{\mathcal{E}}$ and
$g^\flat(\bigoplus_{y \in \mathcal{Y}}y) = \bigoplus_{y \in
\mathcal{Y}}g^\flat(y)$ for each $\mathcal{Y} \subseteq \mathcal{F}$
(i.e., $g^\flat$ is \textit{lower-semicontinuous}).
\end{thm}

\begin{thm}[{\citep[Th. $4.56$]{BICOQ}}]
\label{prop:rescompo} Let  $\mathcal{D}$, $\mathcal{C}$, $\mathcal{B}$ be three semirings. Let  $h : \mathcal{D} \rightarrow \mathcal{C}$
and $f : \mathcal{C} \rightarrow \mathcal{B}$ be residuated
mappings. The following properties hold :
\begin{equation}\label{eq:compResidutriple}
f\circ f^\sharp \circ f =f \text{ and } f^\sharp \circ f \circ f^\sharp=f^\sharp,
\end{equation}
\begin{equation}\label{prop:rescompoFirst}
\begin{array}{lcl}
(f \circ h)^\sharp & = & h^\sharp \circ f^\sharp.
\end{array}
\end{equation}

Let $h : \mathcal{D} \rightarrow \mathcal{C}$
and $g : \mathcal{C} \rightarrow \mathcal{B}$ be dually residuated
mappings. The following properties hold :
\begin{equation}\label{eq:compResidutripledual}
g\circ g^\flat \circ g =g \text{ and } g^\flat \circ g \circ g^\flat=g^\flat,
\end{equation}
\begin{equation}\label{prop:rescompoFirstdual}
\begin{array}{lcl}
(g \circ h)^\flat & = & h^\flat \circ g^\flat.
\end{array}
\end{equation}

\end{thm}

\begin{thm}[{\citep[Th. $4.56$]{BICOQ}}]
\label{prop:rescomposuite} Let  $\mathcal{D}$, $\mathcal{C}$ be two semirings. Let  $h : \mathcal{D} \rightarrow \mathcal{C}$
and $f : \mathcal{D} \rightarrow \mathcal{C}$ be residuated
mappings. The following properties hold :

\begin{equation}\label{eq:inverseorderresiduated}
f \preceq h \Leftrightarrow h^\sharp \preceq f^\sharp,
\end{equation}
\begin{equation}\label{eq:infdistributionresiduated}
(f \oplus h)^\sharp = f^\sharp \wedge h^\sharp.
\end{equation}
Let $h : \mathcal{D} \rightarrow \mathcal{C}$
and $g : \mathcal{D} \rightarrow \mathcal{C}$ be dually residuated
mappings. The following properties hold :
\begin{equation}\label{eq:inverseorderresiduateddual}
g \preceq h \Leftrightarrow h^\flat \preceq g^\flat,
\end{equation}
\begin{equation}\label{eq:infdistributionresiduateddual}
(g \wedge h)^\flat = g^\flat \oplus h^\flat.
\end{equation}
\end{thm}

\begin{thm}[\cite{cohen98a}] \label{thm:imageinclusion}
Let $\sd,\mathcal{C}$ be  semirings, $f: \sd \rightarrow \mathcal{C}$ and $g: \sd \rightarrow \mathcal{C}$ be two residuated mappings, then the following equivalence holds :
$$\im ~f \subset \im g \Leftrightarrow g \circ g^\sharp \circ f= f.$$
\end{thm}
\begin{pf}
If $\im ~f \subset \im ~g $ then  there exists $ \text{ a mapping } h :\mathcal{S} \rightarrow \mathcal{S},   \text { s.t. } f=g \circ h$. According to Equation (\ref{eq:compResidutriple}), $g \circ g^\sharp \circ f=g \circ g^\sharp \circ g \circ h= g \circ h= f$.
If $g \circ g^\sharp \circ f= f$ then $\im f \subset \im ~g$.
\qed
\end{pf}

\begin{proposition}[\cite{cohen98a,cohen97a}, Projection on the image of a mapping]\label{prop:projectiononimage}
Let  $\sd$, $\mathcal{C}$ be  semirings. Let $f: \sd \rightarrow \mathcal{C}$ be a residuated mapping,
mapping $P_f = f\circ f^\sharp$ is a projector and $P_f(c)$ with $c \in \mathcal{C}$ is the greatest element in $\im f$ less than or equal to $c$.
Let $g: \sd \rightarrow \mathcal{C}$ be a dually residuated mapping,
mapping $P_g = g\circ g^\flat$ is a projector and $P_g(d)$ with $d \in \mathcal{C}$ is the lowest element in $\im g$ greater than or equal to $d$.
\end{proposition}
\begin{pf}
According to  Definition \ref{def_residuation}, $P_f(c)=\{ \bigoplus x |f(x) \preceq c \}$ and  $P_g(d)=\{ \bigwedge x | g(x) \succeq d \}$.
According to Equations (\ref{eq:compResidutriple}) and (\ref{eq:compResidutripledual}),
$P_f \circ P_f =f \circ f^\sharp \circ f \circ f^\sharp=f \circ f^\sharp$, and $P_g \circ P_g =g \circ g^\flat \circ g \circ g^\flat=g \circ g^\flat$, hence they are both projectors.
\qed
\end{pf}

The problem of mapping restriction and its connection with
residuation theory is now addressed.
\begin{defn}[Restricted mapping] \label{def_restriction}Let $f:\mathcal{E} \rightarrow \mathcal{F}$ be a
mapping and $\mathcal{A} \subseteq \mathcal{E}$. We will denote
$f_{|\mathcal{A}} : \mathcal{A} \rightarrow \mathcal{F} $ the
mapping defined by $f_{|\mathcal{A}} = f \circ \id_{|\mathcal{A}}$
where $\id_{|\mathcal{A}} : \mathcal{A} \rightarrow \mathcal{E}$
is the canonical injection from $\mathcal{A}$ to $\mathcal{E}$.
Similarly, let $\mathcal{B} \subseteq \mathcal{F}$ with $\im f
\subseteq \mathcal{B}$. Mapping $_{\mathcal{B}|}f : \mathcal{E}
\rightarrow \mathcal{B}$ is defined by $f = \id_{|\mathcal{B}}
\circ {_{\mathcal{B}|}f}$, where $\id_{|\mathcal{B}} : \mathcal{B}
\rightarrow \mathcal{F}$.
\end{defn}

\begin{proposition}[\citep{blyth72}]\label{prop:canonicalinjection}
 Let $\mathcal{S}_{sub}$ be a  complete subsemiring of $\sd$.
Let  $\mathsf{Id}_{|\mathcal{S}_{sub}} : \mathcal{S}_{sub} \rightarrow
 \mathcal{S}$, $x \mapsto x$ be the canonical injection. The injection $\mathsf{Id}_{|\mathcal{S}_{sub}}$
 is both residuated and dually residuated and their residuals are projectors.
\end{proposition}
\begin{pf}
According to Definition \ref{def_continuity}, mapping $\mathsf{Id}_{|\mathcal{S}_{sub}}$ is both $l.s.c.$ and $u.s.c.$, $i.e.$ continuous, and by assumption $\varepsilon \in \mathcal{S}_{sub}$ and $\top \in \mathcal{S}_{sub}$, hence $\mathsf{Id}_{|\mathcal{S}_{sub}}$ is both residuated and dually residuated
(see Theorem \ref{thm:residuation}).
Furthermore, $\mathsf{Id}_{|\mathcal{S}_{sub}} = \mathsf{Id}_{|\mathcal{S}_{sub}} \circ \mathsf{Id}_{|\mathcal{S}_{sub}}$
hence
$(\mathsf{Id}_{|\mathcal{S}_{sub}})^\sharp   = (\mathsf{Id}_{|\mathcal{S}_{sub}} \circ \mathsf{Id}_{|\mathcal{S}_{sub}})^\sharp
= (\mathsf{Id}_{|\mathcal{S}_{sub}})^\sharp \circ (\mathsf{Id}_{|\mathcal{S}_{sub}})^\sharp
$
 which proves that $(\mathsf{Id}_{|\mathcal{S}_{sub}})^\sharp$
 is a projector.
The same can be done for $(\mathsf{Id}_{|\mathcal{S}_{sub}})^\flat$.
\qed
\end{pf}

\begin{proposition} \label{prop:fcodomainedomaine}
Let $f : \mathcal{D} \rightarrow \mathcal{E}$ be a  residuated
mapping, $g : \mathcal{D} \rightarrow \mathcal{E}$ be a  dually residuated
mapping and $\mathcal{D}_{sub}$ (resp. $\mathcal{E}_{sub}$) be a
complete subsemiring of $\mathcal{D}$ (resp. $\mathcal{E}$):
\begin{enumerate}
\item[$1.$] mapping $f_{|\mathcal{D}_{sub}}$ is  residuated and its residual is
given by :
$$ (f_{|\mathcal{D}_{sub}})^\sharp = ( f \circ \mathsf{Id}_{|\mathcal{D}_{sub}})^\sharp = (\mathsf{Id}_{|\mathcal{D}_{sub}})^\sharp \circ
f^\sharp;
$$
\normalsize
\item[$2.$] if $\im f \subset \mathcal{E}_{sub}$ then mapping
$_{\mathcal{E}_{sub}|}f$ is  residuated and its  residual is given
by:
$$
\begin{array}{lcl}
\left ( _{\mathcal{E}_{sub}|} f \right )^\sharp & = & f^\sharp
\circ \mathsf{Id}_{|\mathcal{E}_{sub}} = \left ( f^\sharp \right
)_{|\mathcal{E}_{sub}};
\end{array}
$$
\normalsize
\item[$3.$] mapping $g_{|\mathcal{D}_{sub}}$ is dually residuated and its dual residual is
given by :
$$ (g_{|\mathcal{D}_{sub}})^\flat = ( g \circ \mathsf{Id}_{|\mathcal{D}_{sub}})^\flat = (\mathsf{Id}_{|\mathcal{D}_{sub}})^\flat \circ
g^\flat;
$$
\normalsize
\item[$4.$] if $\im g \subset \mathcal{E}_{sub}$ then mapping
$_{\mathcal{E}_{sub}|}g$ is dually residuated and its dual residual is given
by:
$$
\begin{array}{lcl}
\left ( _{\mathcal{E}_{sub}|} g \right )^\flat & = & g^\flat
\circ \mathsf{Id}_{|\mathcal{E}_{sub}} = \left ( g^\flat \right
)_{|\mathcal{E}_{sub}}.
\end{array}
$$
\normalsize
\end{enumerate}

\end{proposition}
\begin{pf}
Statements 1 and 3 follow directly from Theorem \ref{prop:rescompo} and
Proposition \ref{prop:canonicalinjection}. Statement 2  is obvious
since $f$  is residuated and $\im f \subset \mathcal{E}_{sub}
\subset \mathcal{E}$. Statement 4 can be prove in the same manner.\qed
\end{pf}

\subsection{Closure mappings}

\begin{defn}[Closure mapping]\label{def:closuremapping}
Let $\sd$ be a semiring and $h: \sd \rightarrow \sd$ be an isotone mapping. If $h \circ h=h \succeq \id_\sd$ then
$h$ is a closure mapping. If $h \circ h=h \preceq \id_\sd$ then
$h$ is a dual closure mapping.
\end{defn}
\begin{rem}
According to this definition, it can be checked that the projector $P_f$ (see Proposition \ref{prop:projectiononimage}) is a dual closure mapping, and the projector $P_g$ is a closure mapping.
\end{rem}
\begin{thm}[{\citep[Th. 19 and Th. 20]{cohen98a}}]\label{thm:residuationclosuremapping}
Let $\sd$ be a semiring, $h: \sd \rightarrow \sd$ be a residuated mapping and $g: \sd \rightarrow \sd$ be a dually residuated mapping, then the following equivalences hold:
\small
\begin{equation}\label{eq:closuremappingproperty}
h \text{ is a closure mapping } \Leftrightarrow h^\sharp \text{ is a dual closure mapping} \Leftrightarrow h^\sharp \circ h= h \Leftrightarrow h \circ h^\sharp = h^\sharp ,
\end{equation}
\begin{equation}
g \text{ is a dual closure mapping} \Leftrightarrow g^\flat \text{ is a closure mapping} \Leftrightarrow g \circ g^\flat= g \Leftrightarrow g^\flat \circ g = g^\flat .
\end{equation}
\end{thm}
\normalsize

\begin{proposition}\label{prop:compoclosuremapping}
Let $\sd$ be a semiring, $h: \sd \rightarrow \sd$, $g: \sd \rightarrow \sd$ and $f: \sd \rightarrow \sd$  be three mappings, and assume that $g$ and $f$ are two closure mappings which are residuated. The following equivalence holds
$$ \im h \subset \im f \Leftrightarrow  f \circ h =h, $$
$$g \preceq f \Leftrightarrow f \circ g = f = g^\sharp \circ f \Leftrightarrow \im f \subset \im ~g \Leftrightarrow \im f \subset \im ~g^\sharp .$$
\end{proposition}
\begin{pf}
For the first statement, $ \im h \subset \im f \Rightarrow \exists  ~m \text{ such that }  h= f \circ m  \Rightarrow f\circ h = f \circ f\circ m= f\circ m= h$, since $f$ is a closure mapping, and obviously $f \circ h =h \Rightarrow \im h \subset \im f$.\\
For the second statement, according to the closure mapping definition $g \succeq \id_\sd $, hence $g \succeq \id_\sd \Rightarrow f \circ g \succeq f$. Mapping $f$ is assumed to be a closure mapping, this yields $g \preceq f \Rightarrow  f \circ g \preceq f \circ f=f$. Hence $g \preceq f \Leftrightarrow  f \circ g=f $.

According to Equivalences (\ref{eq:closuremappingproperty}), $g^\sharp$ is a dual closure mapping, therefore  according to Definition \ref{def:closuremapping} $g^\sharp \preceq \id_\sd$, hence $g^\sharp \preceq \id_\sd \Rightarrow g^\sharp  \circ f \preceq f$. According to the assumptions, $f$ and $g$ are residuated, hence Equation (\ref{eq:inverseorderresiduated}) yields $f \succeq g \Leftrightarrow g^\sharp  \succeq f^\sharp \Rightarrow g^\sharp \circ f \succeq f^\sharp \circ f= f$ (the last equality comes from Equivalences (\ref{eq:closuremappingproperty})),
hence $g \preceq f \Leftrightarrow  f =g^\sharp \circ f$.

By considering Equivalences (\ref{eq:closuremappingproperty}) and Theorem \ref{thm:imageinclusion}, $f=g^\sharp \circ f=g\circ g^\sharp \circ f \Rightarrow \im f \subset \im ~g $, on the other hand $\im f \subset \im ~g \Rightarrow \exists ~m \text { such that } f=g\circ m \Rightarrow g^\sharp \circ f =g^\sharp \circ g \circ m = g \circ m= f$.

In the same manner,  Equivalences (\ref{eq:closuremappingproperty}) and Theorem \ref{thm:imageinclusion} yield : $f=g^\sharp \circ f \Rightarrow  \im f \subset \im ~g^\sharp $, on the other hand $\im f \subset \im ~g^\sharp \Rightarrow \exists ~m \text { such that } f=g^\sharp \circ m \Rightarrow  g^\sharp \circ f =g^\sharp \circ g^\sharp \circ m = g^\sharp \circ m= f$ (indeed $g^\sharp=g^\sharp \circ g^\sharp$ since $g^\sharp$ is  a dual closure mapping).
\qed
\end{pf}

\subsection{Applications}


\begin{defn}[Left product, right product]
\label{defn:eq_res} Let $\sd$ be a complete semiring, $a,b \in \sd$,  and $L_a : \sd \rightarrow \sd, x \mapsto a \otimes x$ and $R_a :  \sd \rightarrow \sd, x \mapsto x\otimes a$. Since $\varepsilon$ is absorbing for the multiplicative law and according to distributivity of this law over the additive law, $L_a$ and $R_a$ are both lower semi-continuous, hence both mappings are residuated.  In \cite{BICOQ}, their residuals are
denoted, respectively, by $L^\sharp_a(x)= a \lslash x$ and
$R^\sharp_a(x) = x \rslash a$. Therefore, $a\lslash
b$ (resp. $b \rslash a$) is the greatest solution of $a\otimes x
\preceq b$ (resp. $x \otimes a \preceq b$) and  equality is achieved when $b \in \im L_a$ (resp. $b \in \im R_a$).
It must be noted that $\varepsilon \lslash \varepsilon=\top$ and $\top\lslash \top=\top$.
In the matrix case, mappings $L_A : \sd^{p \times m} \rightarrow \sd^{n\times m}, X \mapsto A \otimes X$ and
$R_A : \sd^{m \times n} \rightarrow \sd^{m \times p}, X \mapsto X\otimes A$
 where $A\in \mathcal{S}^{n \times p}$ , are residuated mappings. The corresponding entries are obtained as follows,
\begin{eqnarray}
 \left ( A \lslash B \right )_{ij} & = & \bigwedge \limits_{k=1
\ldots
n} \left ( a_{ki} \lslash b_{kj} \right ), \label{MatrixLeftResiduation}\\
 \left ( C \rslash A \right )_{ij} & = & \bigwedge \limits_{k=1
\ldots p} \left ( c_{ik} \rslash a_{jk} \right )
\end{eqnarray}
with $B\in \mathcal{S}^{n \times
m}$ and $C\in \mathcal{S}^{m \times p}$.
\end{defn}

\begin{defn}[Kleene star]\label{def:kleenestar}
Let $\sd$ be a complete semiring.
The additive closure of matrix $A \in \sd^{n \times n}$ is defined as follows :
$$K : \sd^{n \times n} \rightarrow \sd^{n \times n}, A \mapsto A^\ast=\bigoplus_{i\in \Nset}A^i, $$
where $A^0=E$, $A^k=A\otimes A^{k-1}$ and  $E$ is the identity matrix, $i.e.$ $\forall i,j \in [1,n]$, $E_{ii}=e$ and $E_{ij}=\varepsilon$ if  $i \neq j$.\\
 This mapping is a closure mapping (indeed $K \circ K=K$ and $K\succeq Id_{S^{n \times n}}$). It is sometimes called the Kleene star operator. Among many references about the Kleene star matrix we can cite \cite{Sergeev09}, where the link between the Kleene star $A^\ast$
 and the \textit{subeigenvectors} of $A$ for an eigenvalue $\lambda$, $i.e.$, vectors $x$  s.t. $A \otimes x \preceq \lambda \otimes x$, was studied.
\end{defn}

\begin{propriete}\label{prop:etoilekleene}
Let $A \in \sd^{n \times n}$, and $X \in \sd^{n \times p}$.
According to  Definition \ref{def:kleenestar} mapping $L_{A^\ast}: \sd^{n \times p} \rightarrow \sd^{n \times p},  X \mapsto A^\ast \otimes X$ is a closure mapping,
(see Definition \ref{def:closuremapping}), hence :
\begin{equation} \label{eq:astarastarastar}
 A^\ast \otimes A^\ast \otimes X=A^\ast \otimes X,
\end{equation}
and as a consequence the following equivalence holds :
\begin{equation} \label{eq:astarxinimastar}
X=A^\ast X \Leftrightarrow X \in \im L_{A^\ast}.
\end{equation}
Furthermore  according to Theorem \ref{thm:residuationclosuremapping},
$L_{A^\ast}^\sharp$ is a dual closure mapping, hence :
\begin{equation} \label{eq:astarsurastar}
A^\ast\lslash A^\ast \lslash X = A^\ast \lslash X,
\end{equation}
according to Equation (\ref{eq:closuremappingproperty}), $L_{A^\ast}\circ L_{A^\ast}^\sharp=L_{A^\ast}^\sharp$ and $L_{A^\ast}^\sharp\circ L_{A^\ast}=L_{A^\ast}$ hence :
 \begin{equation} \label{eq:astartimesastarsurastar}
A^\ast \otimes (A^\ast \lslash X)=A^\ast \lslash X,
\end{equation}
and
\begin{equation} \label{eq:astarsurastartimesastar}
A^\ast \lslash (A^\ast \otimes X)=A^\ast \otimes X.
\end{equation}

According to Proposition \ref{prop:projectiononimage},  Equation (\ref{eq:astartimesastarsurastar}) means that $L_{A^\ast}^\sharp$ is a projector on $\im L_{A^\ast}$.\\
Let $B \in \sd^{n \times n}$ such that $B^\ast \preceq A^\ast$, $i.e.$, $ L_{B^\ast} \preceq L_{A^\ast}$, then according to Proposition \ref{prop:compoclosuremapping}, the following equivalence holds :
\begin{equation} \label{eq:AstarGreaterBstar}
 B^\ast \preceq A^\ast \Leftrightarrow A^\ast  B^\ast  X = A^\ast X =B^\ast \lslash (A^\ast X) \Leftrightarrow \im L_{A^\ast} \subset \im L_{B^\ast} \Leftrightarrow \im L_{A^\ast} \subset \im L_{B^\ast}^\sharp.
 \end{equation}
\end{propriete}
\begin{lem}[\cite{BICOQ}, Lemma 4.77]\label{lem:EquivalenceAstar}
Let $A \in \sd^{n \times n}$, and $X \in \sd^{n \times p}$. The following equivalences hold :
$$X \preceq A \lslash X \Leftrightarrow  X \succeq AX \Leftrightarrow  X=A^\ast X \Leftrightarrow   X=A^\ast\lslash X.$$
\end{lem}
\section{Dual product over semirings} \label{section:dualproduct}
In this section a dual product is considered and its properties are explored.
\begin{defn}[Dual product]\label{def:dualproduct}
 Given a semiring $\sd$, the dual product in $\sd$, denoted $\odot$, is a law assumed to be associative and to have $e$ as neutral element, $i.e.$, $(\sd,\odot,e)$ is a monoid.  Furthermore this dual product is assumed to distribute with respect to $\wedge$ of infinitely many elements, and element  $\top$ is absorbing ($ \forall a,~~\top \odot a = a \odot \top = \top$).
\end{defn}
%


 \begin{defn}[Dual matrix product]\label{def:DualMatrixProduct}
Let $\sd$ be a semiring  and $A\in \mathcal{S}^{n \times p}$, $B\in \mathcal{S}^{p \times m}$ and $C\in \mathcal{S}^{n \times m}$ matrices, then $C=A \odot B$ is defined as :
\begin{eqnarray}\nonumber
C_{ij} =\left ( A \odot B \right )_{ij} & = & \bigwedge \limits_{k=1
\ldots p} \left ( a_{ik} \odot b_{kj} \right ), \label{MatrixDualProduct}
\end{eqnarray}
the identity matrix is denoted $E^\odot$ and is such that $E^\odot_{ii}=e$  and  $E^\odot_{ij}=\top$ for $i\neq j$.\\
\end{defn}
In the sequel, mapping $\Lambda_A :\mathcal{S}^{p \times m} \rightarrow \mathcal{S}^{n \times m}, X \mapsto A\odot X$  will be considered.

\begin{proposition}
Let $\sd$ be a semiring and $A\in \sd^{p\times n}, X\in \sd^{n\times m} $  be matrices, mapping
$\Lambda_A:\sd^{n\times m} \rightarrow \sd^{p\times m},X\mapsto A\odot X$ is upper-semicontinuous, $i.e.$,
\begin{center}
$\Lambda_A(\underset{X\in \mathcal {X}}{\bigwedge}X)=\underset{X\in
\mathcal {X}}{\bigwedge}\Lambda_A(X).$
\end{center}
\end{proposition}

\begin{pf}
Let  $\mathcal {X}$ be a subset of $\sd^{n \times m}$, then according to the definition of $\odot$ the following equalities hold :

$$\begin{array}{cclcl} \Lambda_A(\underset{X\in \mathcal
{X}}{\bigwedge}X) & = &
A\odot(\underset{X\in \mathcal {X}}{\bigwedge}X)\\
(\Lambda_A(\underset{x\in \mathcal {X}}{\bigwedge
}X))_{ij}&=&\underset{k=1}{\overset{n}{\bigwedge }}a_{ik}\odot
(\underset{x\in \mathcal {X}}{\bigwedge } x_{kj})
&=&\underset{k=1}{\overset{n}{\bigwedge }} \underset{x\in
\mathcal {X}}{\bigwedge } (a_{ik}\odot x_{kj})\\
&=&\underset{x\in \mathcal {X}}{\bigwedge }
\underset{k=1}{\overset{n}{\bigwedge }}(a_{ik}\odot x_{kj})
&=&\underset{x\in \mathcal {X}}{\bigwedge }( \Lambda_A(X))_{ij}.
\end{array}$$
\qed
\end{pf}
\begin{cor}\label{cor:dualresiduationmatrixproduct}
Let $\sd$ be a semiring and  $A\in \sd^{n\times p}$ be a matrix. Mapping $\Lambda_A :\mathcal{S}^{p \times m} \rightarrow \mathcal{S}^{n \times m}, X \mapsto A\odot X$ is dually residuated, and its dual residual will be denoted\footnote{This notation was initially introduced
in the talk entitled "Projective $max,+$ semi modules", given by G. Cohen during the International Workshop on $max,+$ Algebra (IWMA Birmingham 2003, in honor of Prof. Cuninghame-Green \cite{cohen03}).}:

\begin{center}
$\Lambda_A^{\flat}:\sd^{n\times m} \rightarrow \sd^{p \times m},X\mapsto A\lslashh X$
\end{center}
with the following rules :
\begin{equation}\label{eq:dualresiduationmatrix}
(A\lslashh
X)_{ij}=\overset{n}{\underset{k=1}\bigoplus}a_{ki}\lslashh
x_{kj},
\end{equation}
and :  $\top\lslashh
x=\varepsilon$, $\varepsilon\lslashh x=\top$ and $\varepsilon\lslashh
\varepsilon=\varepsilon.$
\end{cor}


\begin{proposition}\label{prop:associativitylslashhotimes}
Let $\sd$ be a complete semiring and $A\in \sd^{n\times p}$, $B\in \sd^{n\times
r}$ and $X\in \sd^{p \times q}$ be three matrices. If for each entry $b_{ij}$ of $B$ the following equality holds $b_{ij}\lslashh(a \otimes x)=(b_{ij}\lslashh a) \otimes x,$ $\forall a, x \in \sd$, then the following equality holds :
\begin{equation}
B \lslashh (A \otimes X) =(B\lslashh A) \otimes X.
\end{equation}
\end{proposition}
\begin{pf}
\begin{equation}
\begin{array}{lcl} \nonumber
(B\lslashh(A\otimes
X))_{ij} & = &\overset{n}{\underset{l=1}\bigoplus}b_{li}\lslashh (A\otimes X)_{lj}\\
& = & \overset{n}{\underset{l=1}\bigoplus}b_{li}\lslashh(
\overset{p}{\underset{k=1}\bigoplus}
a_{lk}\otimes x_{kj})\\
& = &\overset{n}{\underset{l=1}\bigoplus}
\overset{p}{\underset{k=1}\bigoplus}
b_{li}\lslashh(a_{lk}\otimes x_{kj}) \text{ since } \Lambda_B^\flat \text{ is lower semi-continuous}\\
&=&
\overset{p}{\underset{k=1}\bigoplus}\overset{n}{\underset{l=1}\bigoplus}(b_{li}\lslashh
a_{lk})\otimes x_{kj} \text{ according to the assumption}  \\
&=&\overset{p}{\underset{k=1}\bigoplus}(B\lslashh
A)_{ik}\otimes x_{kj}=((B\lslashh A)\otimes X)_{ij}.
\end{array}
\end{equation}
\qed
\end{pf}

\begin{defn}\label{def:kleenestarduale}
Let $\sd$ be a semiring. The $\wedge$-closure of $B \in \sd^{n  \times n}$ is defined as:
$$B_\ast=\bigwedge_{k\in \Nset} B^{\odot k},$$
where $B^{\odot 0}=E^\odot$ and $B^{\odot k}=B \odot B^{\odot (k-1)}$.
\end{defn}

\begin{propriete}\label{prop:etoilekleeneduale}
Let $B \in \sd^{n \times n}$, and $X \in \sd^{n \times p}$.
Since $\Lambda_B$ is upper-semicontinuous and, according to  Definition \ref{def:kleenestarduale}, mapping $\Lambda_{B_{\ast}}: \sd^{n \times p} \rightarrow \sd^{n \times p},  X \mapsto B_\ast \odot X$ is a dual closure mapping
(see Definition \ref{def:closuremapping}), hence :
\begin{equation} \label{eq:astarastarastardual}
 B_\ast \odot B_\ast \odot X=B_\ast \odot X,
\end{equation}
and as a consequence the following equivalence holds :
\begin{equation} \label{eq:astarxinimastardual}
X=B_\ast \odot X \Leftrightarrow X \in \im \Lambda_{B_\ast}.
\end{equation}

\end{propriete}


\begin{proposition} \label{propo:equivalenceAdualstar}
Let $\sd$ be a semiring and  $B \in \sd^{n\times n}$ and $X\in \sd^{n\times
p}$ be two matrices. The following statements are equivalent:
\begin{itemize}
 \item[1.] $X\preceq B \odot X$;
 \item[2.] $B\lslashh X\preceq X$;
 \item[3.] $B_\ast \lslashh X= X$;
 \item[4.] $B_\ast \odot X =X $.
\end{itemize}
\end{proposition}

\begin{pf}
$(1)\Rightarrow (2)$ According to Definition \ref{def_residuation} mapping $\Lambda_B^\flat$ is order preserving, hence
$X\preceq B \odot X \Rightarrow B\lslashh X \preceq B\lslashh(B\odot X)$, furthermore the same definition implies $B\lslashh X \preceq B\lslashh(B\odot X)\preceq X$. Hence $X\preceq B \odot X \Rightarrow B\lslashh X\preceq X$.

$(2)\Rightarrow (3)$
According to Equation (\ref{prop:rescompoFirstdual}), $(\Lambda_B \circ \Lambda_B)^\flat=(\Lambda_B)^\flat \circ (\Lambda_B)^\flat$,
hence $B\lslashh(B \lslashh X)=B^{\odot 2}\lslashh X$, furthermore mapping $\Lambda_B^\flat$ is order preserving, then  $$ X \succeq B\lslashh X  \Rightarrow B \lslashh X \succeq B \lslashh(B\lslashh X)= B^{\odot 2}\lslashh X,$$ hence
$$X \succeq B\lslashh X \succeq B^{\odot 2}\lslashh X \succeq \ldots \Rightarrow X \succeq (E^\odot\lslashh X) \oplus (B\lslashh X) \oplus (B^{\odot 2} \lslashh X) \oplus ... $$
furthermore according to Equation (\ref{eq:infdistributionresiduateddual}) and to Definition \ref{def:kleenestarduale}, $$  (E^\odot\lslashh X) \oplus (B\lslashh X) \oplus (B^{\odot 2} \lslashh X) \oplus ...= (E^\odot \wedge B \wedge B^{\odot 2}\wedge ... ) \lslashh X=B_\ast \lslashh X,$$ then, $ X \succeq (B\lslashh X) \Rightarrow  X \succeq B_\ast \lslashh X.$ On the other hand $B_\ast \preceq E^\odot$ then
$B_\ast \lslashh X \succeq X$, hence $ X \succeq (B\lslashh X) \Rightarrow  X =B_\ast \lslashh X.$

$(3)\Rightarrow (4)$ From Definition \ref{def_residuation} (Equation (\ref{dualresiduated1})) the following inequality holds : $ B_\ast \odot (B_\ast \lslashh X) \succeq X $, hence,
 $$X=B_\ast \lslashh X \Rightarrow B_\ast \odot X=B_\ast \odot (B_\ast \lslashh X)\succeq X,$$ but the definition of the dual closure
yields $B_\ast \odot X \preceq X$, hence   $$X=B_\ast \lslashh X \Rightarrow B_\ast \odot X = X.$$

$(4)\Rightarrow (1)$ According to Definitions \ref{def:dualproduct} and \ref{def:kleenestarduale}, Mapping $\Lambda_{B_\ast}$ is upper semicontinuous, then
$$ B_\ast \odot X =( E^\odot \wedge B \wedge B^{\odot 2}\wedge ...)\odot X= ( X  \wedge B\odot X \wedge B^{\odot 2} \odot X\wedge ...),$$
hence $X=B_\ast \odot X \Rightarrow X \preceq B \odot X.$
\qed
\end{pf}
\section{The Inequality $ A\otimes X \preceq X \preceq B \odot X $}\label{section:axbx}
\begin{proposition}
Let $\sd$ be a semiring and $A, B  \in \sd^{n \times n}$ and $X \in \sd^{n \times m}$.
The following equivalence holds :
\begin{equation}
A\otimes X \preceq X \preceq B \odot X \Leftrightarrow X \in \im L_{A^\ast} \cap \im \Lambda_{B_\ast}.
\end{equation}
\end{proposition}
\begin{pf}
Direct application of Equivalence (\ref{eq:astarxinimastar}) (see Property \ref{prop:etoilekleene}) and of  Equivalence (\ref{eq:astarxinimastardual}) (see Property \ref{prop:etoilekleeneduale}).
\qed
\end{pf}

\begin{proposition}\label{propo:ProjectorInterBstarAstar}
Let $\sd$ be a semiring and $A, B  \in \sd^{n \times n}$ and $X \in \sd^{n \times m}$.\\ If $\forall X$,  the  equality  $B_\ast \lslashh (A^\ast \otimes X)=(B_\ast \lslashh A^\ast) \otimes X$ holds, then the mapping
$$P : \sd^{n \times m} \rightarrow \sd^{n \times m}, X \mapsto (B_\ast \lslashh A^\ast)^\ast \lslash X,$$
is a projector in $\im L_{A^\ast} \cap \im \Lambda_{B_\ast}$, formally  $$P(X)=\{ \bigvee Y | Y\preceq X  \text{ and } Y                                                                                                                                                                         \in \im L_{A^\ast} \cap \im \Lambda_{B_\ast}  \}.$$
\end{proposition}
\begin{pf}

First, according to Equations (\ref{eq:astarsurastar}) and (\ref{eq:astartimesastarsurastar}), $P$ is a projector on the image of $L_{(B_\ast \lslashh A^\ast)^\ast}$, and $P(X)\preceq X$.\\
According to  Definition \ref{def:kleenestarduale}, $B_\ast \preceq E^\odot$, then  $B_\ast \lslashh A^\ast \succeq E^\odot \lslashh A^\ast=A^\ast$ and $(B_\ast \lslashh A^\ast)^\ast \succeq (A^\ast)^\ast=A^\ast$, which, according to Equation (\ref{eq:AstarGreaterBstar}), implies that $\im L_{(B_\ast \lslashh A^\ast)^\ast} \subset \im L_{A^\ast}$, hence $P(X) \in \im L_{A^\ast}$.\\
Since $P(X) \in \im L_{(B_\ast \lslashh A^\ast)^\ast}$, equality $P(X)=(B_\ast \lslashh A^\ast)^\ast P(X)$ holds, and according to Lemma \ref{lem:EquivalenceAstar}, this is equivalent to
$ P(X) \succeq (B_\ast \lslashh A^\ast) \otimes P(X).$\\
Because of the assumption, the  equality :\\
$(B_\ast \lslashh A^\ast) \otimes P(X)=B_\ast \lslashh (A^\ast \otimes P(X))$ holds, furthermore $P(X) \in \im L_{A^\ast}$, therefore $A^\ast \otimes P(X)=P(X)$, hence  $$P(X) \succeq (B_\ast \lslashh A^\ast) \otimes P(X)=B_\ast \lslashh (A^\ast \otimes P(X))=B_\ast \lslashh  P(X).$$
Otherwise, $B_\ast \preceq E^\odot$, then $$B_\ast \lslashh  P(X) \succeq E^\odot \lslashh P(X) =P(X).$$
Hence, $ P(X) = B_\ast \lslashh  P(X)$.\\
Furthermore, Proposition \ref{propo:equivalenceAdualstar} gives  :
$$P(X) = B_\ast \lslashh  P(X) = B_\ast \odot  P(X),$$
then, by considering Equivalence (\ref{eq:astarxinimastardual}), this implies that $P(X) \in \im \Lambda_{B_\ast}.$

Now we show that $P(X)$ is the greatest element in  $\im L_{A^\ast} \cap \im \Lambda_{B_\ast}$ less or equal to $X$.

Let  $Y \in  \im L_{A^\ast} \cap \im \Lambda_{B_\ast}$ such that $Y \preceq X$, hence according to Lemma   \ref{lem:EquivalenceAstar} and Proposition \ref{propo:equivalenceAdualstar}, the following equalities hold :
$$Y=A^\ast \otimes Y=B_\ast \odot Y= B_\ast \lslashh Y= B_\ast \lslashh (A^\ast Y),$$ and because of the assumption $B_\ast \lslashh (A^\ast Y) =(B_\ast \lslashh A^\ast) Y$.\\
From Definition \ref{def_residuation}, $Y=(B_\ast \lslashh A^\ast) Y \Rightarrow Y \preceq (B_\ast \lslashh A^\ast) \lslash Y $, and from Lemma \ref{lem:EquivalenceAstar}, this is equivalent to $Y=(B_\ast \lslashh A^\ast)^\ast \lslash Y$. Mapping $L_{(B_\ast \lslashh A^\ast)^\ast}^\sharp$ being an isotone mapping, the following implication holds : $ Y \preceq X \Rightarrow (B_\ast \lslashh A^\ast)^\ast \lslash Y \preceq (B_\ast \lslashh A^\ast)^\ast \lslash X$ which means that if $Y \preceq X$ then $Y=(B_\ast \lslashh A^\ast)^\ast \lslash Y\preceq P(X)$.\qed
\end{pf}
\begin{rem}
The previous result shows that $P(X_0)$ is the greatest solution of the following system of inequalities
$$ A\otimes X \preceq X \preceq B \odot X \text{ and } X \preceq X_0,$$
which is equivalent to
$$ A^\ast \otimes X=B_\ast \odot X=X \text{ and } X\preceq X_0.$$
This projector can be useful to synthesize a controller for manufacturing systems subject to constraints. This kind of problem  can be seen as a model matching problem (see \cite{Shang09,Shang11}) and is of practical interest in many industrial applications (see e.g. \cite{brunsch11} for an example from high-throughput-screening).
\end{rem}
\section{Examples}\label{section:Example}
The results introduced in the previous section are illustrated in two semirings of practical interest in control theory of discrete
event systems.
\begin{defn}[Semiring $\zmax$]
According to Definition \ref{def:semiring}, the set $\overline{\Zset} = \Zset \cup \{ - \infty, + \infty \}$
endowed with the max operator as $\oplus$ and the classical sum
 as $\otimes$ is a complete idempotent semiring, denoted $\zmax$, where $\varepsilon = - \infty$,
$e=0$ and $\top = + \infty$. The greatest lower bound is $a \wedge
b=min(a,b)$, and $b\lslash a=a-b$. Furthermore  $a\odot b=a+b$ and $b\lslashh a=a-b$.
Hence, except $\varepsilon$ and $\top$,  all elements  admit a multiplicative inverse
$a^{-1}$, $i.e$, $a\otimes a^{-1}=a^{-1}\otimes a=e$ and $a\odot a^{-1}=a^{-1}\odot a=e$.
As a consequence, the following distributivity properties hold : $c\otimes (a\wedge b)=(c\otimes a) \wedge (c \otimes b)$, $c \odot (a\wedge b)=(c\odot a) \wedge (c \odot b)$ and $c \odot (a\oplus b)=(c\odot a) \oplus (c \odot b)$.
Obviously, this is not true in the matrix case.
\end{defn}
\begin{exmp}
Let $A=\begin{pmatrix}1 & \top & 3 \\ 4 & \varepsilon & 6\end{pmatrix}$ and
$B=\begin{pmatrix}8 \\ 9 \\ 10\end{pmatrix}$, $C=\begin{pmatrix}1 & 2 \\ 3& 4 \\ 5 & 6\end{pmatrix}$ be matrices with
entries in $\zmax$. The  product of these matrices is:
$$ A \otimes B =\begin{pmatrix}
(1\otimes 8) \oplus (\top \otimes 9) \oplus (3 \otimes 10) \\
(4\otimes 8) \oplus (\varepsilon \otimes 9) \oplus (6 \otimes 10)
\end{pmatrix} = \begin{pmatrix}\top \\ 16 \end{pmatrix},$$
and the dual product yields
$$ A \odot B =\begin{pmatrix}
(1\odot 8) \wedge (\top \odot 9) \wedge (3 \odot 10) \\
(4\odot 8) \wedge (\varepsilon \odot 9) \wedge (6 \odot 10)
\end{pmatrix} = \begin{pmatrix}9 \\ \varepsilon  \end{pmatrix}.$$
The greatest  solution of $C\otimes X \preceq B$ is given by
$$ C \lslash B =\begin{pmatrix}
(1\lslash 8) \wedge (3 \lslash 9) \wedge (5 \lslash 10) \\
(2\lslash 8) \wedge (4 \lslash 9) \wedge (6 \lslash 10)
\end{pmatrix} = \begin{pmatrix}5 \\ 4 \end{pmatrix},$$
and the smallest  solution of $C\odot X \succeq B$ is given by
$$ C \lslashh B =\begin{pmatrix}
(1\lslashh 8) \oplus (3 \lslashh 9) \oplus (5 \lslashh 10) \\
(2\lslashh 8) \oplus (4 \lslashh 9) \oplus (6 \lslashh 10)
\end{pmatrix} = \begin{pmatrix}7 \\ 6 \end{pmatrix}.$$
\end{exmp}
\begin{rem}
The dual product can be used to perform residuation of matrices in the $(max,plus)$ algebra (see \cite{cuning80}).
More precisely, in this particular case, $A\lslash B=-A^T \odot B$.
\end{rem}

\begin{defn}[Semiring $\zmaxgamma$, \cite{BICOQ}, $\S 5.3.2$]\label{def:zmaxgamma}
According to Definition \ref{def:formalpowerseries}, the set of  non-decreasing formal power series in one variable
$\gamma$ with coefficients in the semiring $\zmax$ and exponents in $\Zset$ is
a semiring denoted $\zmaxgamma$, where $\gamma^\ast=\bigoplus_{ i \in \Nset} \gamma^i$ (see Definition \ref{def:kleenestar}).
The neutral element of addition is the series $\varepsilon(\gamma)=\bigoplus_{ k \in \Zset} \varepsilon \gamma^k$
and the neutral element of multiplication is the series $e(\gamma)=\bigoplus_{ k \in \Nset} e \gamma^k$, furthermore
$\top(\gamma)=\bigoplus_{ k \in \Zset} \top \gamma^k$.  The monomials are defined as $\gamma^\ast(t\gamma^n)=\bigoplus_{ k \in \Nset} t \gamma^{n+k}$. In order to keep notation simple, this will be denoted $t\gamma^n$ in the sequel of this paper. In the same way, a series will be simply denoted $s=\bigoplus_{i \in I_S} t_i \gamma^{n_i}$, where $I_S\subset \Nset$.
The computational rules between monomials are the following :
\small
 \begin{eqnarray}
t_1 \gamma^n \oplus t_2 \gamma^n=max(t_1,t_2)\gamma^n, & t_1 \gamma^{n_1} \otimes t_2 \gamma^{n_2}=(t_1+t_2)\gamma^{n_1+n_2},\label{eq:oplusotimemonomials}\\
t_1 \gamma^{n_1} \wedge t_2 \gamma^{n_2}=min(t_1,t_2)\gamma^{max(n_1,n_2)}, &  t_1 \gamma^{n_1} \odot t_2 \gamma^{n_2}=(t_1+t_2)\gamma^{n_1+n_2},\label{eq:wedgeodotmonomials}\\
(t_1 \gamma^{n_1}) \lslash (t_2 \gamma^{n_2})=(t_2-t_1)\gamma^{n_2-n_1},  & (t_1 \gamma^{n_1}) \lslashh (t_2 \gamma^{n_2})=(t_2-t_1)\gamma^{n_2-n_1}\label{eq:lslashlslahhmonomials}.
\end{eqnarray}
\normalsize
Furthermore, the order relation is such that $t_1\gamma^{n_1}\succeq t_2 \gamma^{n_2}\Leftrightarrow n_1 \leq n_2 \text{ and } t_1 \geq t_2$.
 According to these rules, a non decreasing series admits many representations ($e.g.$, $2\gamma^2\oplus 3\gamma^2=3\gamma^2$) and one of which is canonical. It is the representation whose $t_0<t_1<...$ and $n_0<n_1<...$.
The computation rules between two series $s=\bigoplus_{i \in I_S} t_i \gamma^{n_i} $ and $s'= \bigoplus_{j \in I_{S'}} t_j \gamma^{n_j}$ are  given as follows :
\begin{eqnarray}\label{operationseries}
  s \oplus s' &= &\bigoplus_{i \in I_S} t_i \gamma^{n_i} \oplus \bigoplus_{j \in I_{S'}} t_j \gamma^{n_j}, \label{eq:sommeseries}\\
  s \otimes s' &= & \bigoplus_{i\in I_S}\bigoplus_{j\in I_{S'}} (t_i+t_j)\gamma^{n_i+n_j}, \label{eq:productseries}\\
  s \wedge s' &= & \bigoplus_{i\in I_S}\bigoplus_{j\in I_{S'}} min(t_i,t_j)\gamma^{max(n_i,n_j)}, \label{eq:wedgeseries}\\
  s \rslash s'=  s' \lslash s &= & \bigwedge_{j\in I_{S'}}\bigoplus_{i\in I_S} (t_i-t_j)\gamma^{n_i-n_j}\label{eq:residuationseries}.
\end{eqnarray}
According to Definition \ref{def:dualproduct}, the dual product has to distribute with respect to the operator $\wedge$,
hence it is only defined between a monomial and a series in the following way :
\begin{eqnarray}\label{eq:operationseriesodot}
 t\gamma^n \odot s= \bigoplus_{i \in I_S}(t+t_i)\gamma^{n+n_i}.
\end{eqnarray}
It can be checked that $a\odot(s\wedge s')=(a\odot s)\wedge (a \odot s')$. The dual residual is then given by :
\begin{eqnarray}\label{eq:operationserieslslashh}
 t\gamma^n \lslashh s= \bigoplus_{i \in I_S}(t_i-t)\gamma^{n_i-n}.
\end{eqnarray}
In \cite{cohen89a}, periodic series were introduced. They are defined as  $s=p \oplus q \otimes r^\ast$ where
$p=\overset{m}{\underset{i=1}\bigoplus}t_i \gamma^{n_i}$ (respectively $q=\overset{l}{\underset{j=1}\bigoplus}t_j \gamma^{n_j}$) is a polynomial depicting the transient (resp. the periodic) behavior, and
$r=\tau\gamma^\nu$ is a monomial depicting the periodicity
allowing to define the asymptotic slope of the series as
$\sigma_{\infty}(s)=\nu/\tau$. Sum, product, Kleene star and residuation of periodic series are periodic series
(see \cite{gaubert92a}), and algorithms and software toolboxes are
available in order to handle them (see \cite{lhommeau00a}). In the same way, the dual product and its dual residual are well defined.
Below, only properties concerning asymptotic slopes are recalled:
\begin{equation}
\begin{array}{lll}
\nonumber \sigma_{\infty}(s \oplus s')& = &\min(\sigma_{\infty}(s), \sigma_{\infty}(s')),\\
\nonumber \sigma_{\infty}(s \otimes  s') & = &\min(\sigma_{\infty}(s),\sigma_{\infty}(s')),\\
\nonumber \sigma_{\infty}(s \wedge  s') & = & \max(\sigma_{\infty}(s),\sigma_{\infty}(s')), \\
\nonumber \sigma_{\infty}(s^\ast) & = & \min(
\underset{i=1..m} {min}(n_i/t_i),\underset{j=1..l}{min}(n_j/t_j),\sigma_{\infty}(s)),\\
\sigma_{\infty}(t \gamma^n \odot s)& = & \sigma_{\infty}(s)\\
\sigma_{\infty}((t\gamma^n) \lslashh s)& = & \sigma_{\infty}(s)\\
\nonumber \textnormal{if~} \sigma_{\infty}(s) \leq \sigma_{\infty}(s')&
\textnormal{~then~}& \sigma_{\infty}(s' \lslash s) =
\sigma_{\infty}(s), \textnormal{~else~~} s' \lslash s=\varepsilon.
\label{penteresiduation}\\
\end{array}
\end{equation}

\end{defn}
\begin{exmp}\label{exmp:B_astmonomials}
Let $B=\begin{pmatrix}\top & 15\gamma^3 & 7 \gamma^0 & \top \\
\top & \top & \top & \top \\
3\gamma^0 & 8\gamma^4 & \top &  \top \\
6\gamma & 4 \gamma^5 & \top & \top\end{pmatrix}$ be a matrix where the entries are monomials in $\zmaxgamma$.
According to Definitions \ref{def:zmaxgamma} and \ref{def:DualMatrixProduct}. It can be checked that :\\
$B^{\odot 2}=\begin{pmatrix}10\gamma^0 & 15\gamma^4 & \top  & \top \\
\top & \top & \top & \top \\
\top & 18\gamma^3 & 10\gamma^0 &  \top \\
\top & 21 \gamma^4 & 13\gamma & \top\end{pmatrix}$ and $B^{\odot 3}=\begin{pmatrix}\top & 25\gamma^3 & 17\gamma^0 & \top \\
\top & \top & \top & \top \\
13\gamma^0 & 18\gamma^4 & \top &  \top \\
16\gamma & 21 \gamma^5 & \top & \top\end{pmatrix}$ \\
It can be also checked that $B^{\odot n}\succeq B^{\odot 3} \forall n>3$, hence :\\
 $B_\ast=E^\odot \wedge B \wedge B^{\odot 2} \wedge B^{\odot 3}\wedge \ldots =\begin{pmatrix}e & 15\gamma^4 & 7\gamma^0 & \top \\
\top & e & \top & \top \\
3\gamma^0 & 8\gamma^4 & e &  \top \\
6\gamma & 4 \gamma^5 & 13\gamma & e\end{pmatrix}$.\\
Note that, due to the computation rules (\ref{eq:wedgeodotmonomials}), the entries of matrix $B_\ast$ are always monomials.
\end{exmp}

\begin{rem}
From these examples, it can be seen that the assumption of Proposition \ref{prop:associativitylslashhotimes},
$i.e.$,  that $b_{ij}\lslashh(a \otimes x)=(b_{ij}\lslashh a) \otimes x$, is clearly fulfilled in the semiring $\zmax$ (indeed $b_{ij}-(a+x)=(b_{ij}-a)+x$).
In $\zmaxgamma$, the dual product is only defined between monomials and series. Hence by considering monomial $b_{ij}=t\gamma^n$, series
$a=\bigoplus_{i\in I_A}t_i\gamma^{n_i}$ and $x=\bigoplus_{j\in I_X}t_j\gamma^{n_j}$, and
according to Equations (\ref{eq:productseries}) and (\ref{eq:operationserieslslashh}) the following equalities hold :
\begin{eqnarray}
\nonumber (t\gamma^n)\lslashh(a \otimes x)&=&(t\gamma^n)\lslashh(\bigoplus_{i\in I_A} t_i\gamma^{n_i}\otimes \bigoplus_{j\in I_X}t_j\gamma^{n_j}) =(t\gamma^n)\lslashh(\bigoplus_{i\in I_A}\bigoplus_{j\in I_X} (t_i+t_j)\gamma^{n_i+n_j})\\
\nonumber &=&\bigoplus_{i\in I_A}\bigoplus_{j\in I_X} (t_i+t_j-t)\gamma^{n_i+n_j-n}=(\bigoplus_{i\in I_A} (t_i-t)\gamma^{n_i-n})\otimes\bigoplus_{j\in I_X}t_j\gamma^{n_j}\\
\nonumber &=&((t\gamma^n)\lslashh a)\otimes x.
\end{eqnarray}

The assumption $B_\ast\lslashh (A^\ast \otimes X)=(B_\ast\lslashh A^\ast) \otimes X$ used in Proposition \ref{propo:ProjectorInterBstarAstar} is still valid in $\zmax$  since $B_\ast$ is with entries in the semiring  $\zmax$. In the same way, it also holds in $\zmaxgamma$ since all entries of $B$  are assumed to be monomials and, as noticed in Example \ref{exmp:B_astmonomials}, under this assumption all entries of $B_\ast$ are monomials.
\end{rem}
\section{Interval Analysis over idempotent semirings}\label{section:intervalanalysis}
Interval mathematics was pioneered by R.E. Moore (see \cite{moore79}) as a tool for
bounding rounding errors in computer programs. Since then, interval mathematics has been developed into a general methodology for investigating numerical uncertainty in many problems and algorithms \cite{Jaulinbook01}.
In \cite{litvinov01c} idempotent semirings were extended to interval arithmetic (see also \cite{Myskova12}).
Below some preliminary statements are recalled from this reference.
\begin{defn}[Interval] \label{def:closedinterval}
Let $\sd$ be a semiring. A (closed) interval is a set of the form
$\mathbf{x} = [\underline{x},\overline{x}] = \{t \in \mathcal{S} |
\underline{x} \preceq t \preceq \overline{x} \}$, where
 $\underline{x} \in \sd$ and  $\overline{x} \in \sd$ (with $\underline{x}\preceq \overline{x}$) are called the lower
and the upper bounds of the interval $\mathbf{x}$, respectively.
\end{defn}
\begin{defn}[Semiring of intervals]
 The set of intervals denoted by $\Is$,
endowed with the following element-wise algebraic operations
\begin{equation}\label{eq:SumProductInterval}
\begin{array}{lcl}
\mathbf{x} \stackrel{-}{\oplus} \mathbf{y} \triangleq  \left
[\underline{x} \oplus \underline{y}, \overline{x} \oplus
\overline{y} \right ] & \textnormal{~~~and~~~} & \mathbf{x}
\stackrel{-}{\otimes} \mathbf{y} \triangleq \left [\underline{x}
\otimes \underline{y}, \overline{x} \otimes \overline{y} \right ]
\end{array}
\end{equation}
is a semiring, where the intervals
$\pmb{\varepsilon}=[\varepsilon,\varepsilon]$ and
$\mathbf{e}=[e,e]$ are the neutral elements of
$\Is$. The canonical order $\preceq_{\Is}$ induced by the additive law is such that
$ \mathbf{x} \stackrel{-}{\oplus} \mathbf{y}= [\underline{x} \oplus \underline{y}, \overline{x} \oplus
\overline{y} ]\Leftrightarrow \mathbf{x} \preceq_{\Is} \mathbf{y} \Leftrightarrow \underline{x} \preceq_\sd \underline{y}$ and $\overline{x} \preceq_\sd \overline{y}$, where $\preceq_\sd$ is the order relation in $\sd$.
\end{defn}
\begin{rem}
In the sequel, in the absence of ambiguity, the order relation in $\Is$ will be  denoted $\preceq$.
Operations (\ref{eq:SumProductInterval}) give the tightest intervals containing all results of the same operations to arbitrary elements of its interval operands.
\end{rem}
\begin{rem} \label{def:sommeinfinie}
 Let $\mathcal{S}$ be a complete semiring and $\{\mathbf{x}_{\alpha}\}$ be an infinite subset of
$\Is$, the infinite sum of elements of this
subset is :
$$ \overline{\bigoplus_{\alpha}} \mathbf{x}_{\alpha} = \left [ \bigoplus_{\alpha} \underline{x}_{\alpha}, \bigoplus_{\alpha} \overline{x}_{\alpha} \right
].
$$
The top element
is given by
$\pmb{\top}=[\top,\top]$.
\end{rem}
\begin{rem}
Note that if $\mathbf{x}$ and $\mathbf{y}$ are intervals in
$\Is$, then $\mathbf{x} \subset \mathbf{y}$
iff $\underline{y} \preceq \underline{x} \preceq \overline{x}
\preceq \overline{y}$. In particular, $\mathbf{x} = \mathbf{y}$
iff $\underline{x} =
\underline{y}$ and $\overline{x} = \overline{y}$.
\end{rem}
\begin{rem}
An interval for which $\underline{x} = \overline{x}$ is called
\textit{degenerate}. Degenerate intervals allow to represent
numbers without uncertainty. In this case $\mathbf{x}$ will be simply denoted $x$.
\end{rem}
\begin{rem}
$\Is$ is not a semifield even if $\sd$ is one. Indeed, except for degenerate intervals, an interval does not admit a multiplicative inverse.
\end{rem}

\begin{defn}[Dual product over semiring $\Is$]
 In a semiring of intervals, the dual product $\odot$ is defined as :
 $ \mathbf{x} \stackrel{-}{\odot} \mathbf{y} \triangleq  \left
[\underline{x} \odot \underline{y}, \overline{x} \odot
\overline{y} \right ]$, where $\odot$ is the dual product in $\sd$.

 \end{defn}

In \cite{Litvinov02} (see also \cite{lhommeau04}), it has been shown that order preserving mappings admit a natural extension over the semirings of intervals by considering the image of the interval bounds in an independent way. Especially the additive closure and $\wedge$-closure can be computed in an efficient way and are defined as follows.
\begin{proposition}[\cite{Litvinov02},\cite{lhommeau04}] \label{prop:intervalclosuredualclosure}
Let $\Is$ be a semiring of intervals. The additive closure of matrix $\mathbf{A} \in \Is^{n \times n}$ is given by :
$$\mathbf{A}^\ast=[ \underline{A},\overline{A} ]^\ast= [ \underline{A}^\ast,\overline{A}^\ast],$$
and its $\wedge$-closure is :
$$\mathbf{A}_\ast=[ \underline{A},\overline{A} ]_\ast= [ \underline{A}_\ast,\overline{A}_\ast].$$
\end{proposition}

\begin{notation}[Semiring of pairs]
Let $\mathcal{S}$ be a complete semiring. The set of pairs $(x',x'')$ with $x' \in
\mathcal{S}$ and $x'' \in \mathcal{S}$ is a complete semiring denoted by $\cd$ with
$(\varepsilon,\varepsilon)$ as the zero element, $(e,e)$ as the
identity element and $(\top,\top)$ as top element (see Definition \ref{def:semiring}).
The set of pairs  $(x',x'')$ such that $x' \preceq x''$ is a complete subsemiring of $\cd$ (see Definition \ref{def:subsemiring}).
It will be denoted $\cod$.
\end{notation}

\begin{proposition}\label{prop:reisudationinjection}
The canonical injection
$\mathsf{Id}_{|\mathcal{C}_\mathtt{O}(\mathcal{S})} :
\mathcal{C}_\mathtt{O}(\mathcal{S}) \rightarrow \cd$ is both residuated and dually
residuated. Its  residual
$(\mathsf{Id}_{|\mathcal{C}_\mathtt{O}(\mathcal{S})})^\sharp$ is a
projector. Its practical
computation is given by :
\begin{equation}\label{eq:residuationpair}
\begin{array}{lcl}
(\mathsf{Id}_{|\mathcal{C}_\mathtt{O}(\mathcal{S})})^\sharp((x',x'')) & = &
(x'\wedge x'' ,  x'') = (\widetilde{x}', \widetilde{x}'').
\end{array}
\end{equation}
Its dual residual
$(\mathsf{Id}_{|\mathcal{C}_\mathtt{O}(\mathcal{S})})^\flat$ is a
projector. Its practical
computation is given by :
\begin{equation}\label{eq:dualresiduationpair}
\begin{array}{lcl}
(\mathsf{Id}_{|\mathcal{C}_\mathtt{O}(\mathcal{S})})^\flat((x',x'')) & = &
(x' , x' \oplus x'') = (\widetilde{x}', \widetilde{x}'').
\end{array}
\end{equation}

\normalsize
\end{proposition}
\begin{pf}
This theorem is a direct application of Proposition
\ref{prop:canonicalinjection}, since
$\mathcal{C}_\mathtt{O}(\mathcal{S})$ is a subsemiring of $\cd$.
Practically, let us consider $(x',x'') \in \cd$, we have
$(\mathsf{Id}_{|\mathcal{C}_\mathtt{O}(\mathcal{S})})^\sharp((x',x'')) =
(\widetilde{x}', \widetilde{x}'') = (x'\wedge x'' ,  x'')$, which
is the greatest pair such that :

$$   \widetilde{x}' \preceq x', ~~~\widetilde{x}'' \preceq x'' ~~\textnormal{and}~~\widetilde{x}' \preceq \widetilde{x}''.$$
On the other hand, we have $(\mathsf{Id}_{|\mathcal{C}_\mathtt{O}(\mathcal{S})})^\flat((x',x'')) =
(\widetilde{x}', \widetilde{x}'') = (x' ,x' \oplus  x'')$, which
is the smallest pair such that :

$$   \widetilde{x}' \succeq x', ~~~\widetilde{x}'' \succeq x'' ~~\textnormal{and}~~\widetilde{x}'' \succeq \widetilde{x}'.$$\qed
\end{pf}

\begin{proposition}[\cite{hardouin09b}] \label{prop:residuation_couple_ordonnes}
Mapping $L_{(a',a'')} : \cod \rightarrow
\cod,(x',x'') \mapsto  (a'\otimes x',a''\otimes x'')$ with $(a',a'') \in \cod$ is residuated. Its residual is
equal to
\begin{equation}\label{eq:residuationcouplesorodonnes}
L_{(a',a'')}^\sharp : \cod
\rightarrow \cod, (x',x'') \mapsto
(a' \lslash x' \wedge a'' \lslash x'', a'' \lslash x'').
\end{equation}
\end{proposition}

\begin{proposition}[\cite{hardouin09b}] \label{res_interval}
Let $\Is$ be a semiring of intervals. Mapping $L_{\mathbf{a}} : \Is \rightarrow
\Is,\mathbf{x} \mapsto \mathbf{a}
\stackrel{-}{\otimes} \mathbf{x}$ is residuated. Its residual is
equal to $$L_{\mathbf{a}}^\sharp : \Is
\rightarrow \Is, \mathbf{x} \mapsto
\mathbf{a} \overline{\lslash} \mathbf{x}= [\underline{a}\lslash
\underline{x} \wedge \overline{a}\lslash\overline{x} ,
\overline{a} \lslash \overline{x}].$$ Therefore, $ \mathbf{a}
\overline{\lslash} \mathbf{b}$ is the greatest solution of
$\mathbf{a} \stackrel{-}{\otimes} \mathbf{x} \preceq \mathbf{b} $,
and equality is achieved if $\mathbf{b} \in \im L_{\mathbf{a}}$.
\end{proposition}
\begin{rem}
In the same manner, it can be shown that mapping $R_\mathbf{a} :
\Is \rightarrow \Is,
\mathbf{x} \mapsto \mathbf{x} \stackrel{-}{\otimes} \mathbf{a}$ is residuated.
\end{rem}

\begin{proposition} \label{prop:dualresiduation_couple_ordonnes}
Mapping $\Lambda_{(a',a'')} : \cod \rightarrow
\cod, (x',x'') \mapsto  (a'\odot x',a''\odot x'')$ with $(a',a'') \in \cod$ is dually residuated. Its dual residual is
equal to
\begin{equation}\label{eq:dualresiduationcouplesorodonnes}
\Lambda_{(a',a'')}^\flat : \cod
\rightarrow \cod, (x',x'') \mapsto
(a' \lslashh x', a' \lslashh x' \oplus a'' \lslashh x'').
\end{equation}
\end{proposition}
\begin{pf}
According to Corollary \ref{cor:dualresiduationmatrixproduct}, mapping $\Lambda_{(a',a'')} : \cd \rightarrow
\cd,(x',x'') \mapsto  (a'\odot x',a''\odot x'')$ is dually residuated and its dual residual is
$\Lambda_{(a',a'')}^\flat : \cd
\rightarrow \cd, (x',x'') \mapsto
(a' \lslashh x' , a'' \lslashh x'')$.  Mapping  $\Lambda_{(a',a'')}$ is order preserving, hence $\im \Lambda_{(a',a'')|\cod} \in \cod$. Furthermore, the canonical injection
$\mathsf{Id}_{|\cod} :
\cod \rightarrow \cd$ is dually
residuated. Hence Proposition \ref{prop:fcodomainedomaine} yields $$(_{ \cod |}\Lambda_{(a',a'') |\cod})^\flat=(_{ \cod |}\Lambda_{(a',a'')} \circ  \mathsf{Id}_{|\mathcal{C}_\mathtt{O}(\mathcal{S})})^\flat=(\mathsf{Id}_{| \cod})^\flat \circ (\Lambda_{(a',a'')})^\flat \circ \mathsf{Id}_{|\mathcal{C}_\mathtt{O}(\mathcal{S})}.$$
To conclude, Equation (\ref{eq:dualresiduationpair}) of Proposition \ref{prop:reisudationinjection} yields equation (\ref{eq:dualresiduationcouplesorodonnes}).
\qed
\end{pf}

\begin{proposition} \label{dualres_interval}
Let $\sd$ be a semiring  and $\Is$ be a semiring of intervals. Mapping $\Lambda_{\mathbf{a}} : \Is \rightarrow
\Is,\mathbf{x} \mapsto \mathbf{a}
\stackrel{-}{\odot} \mathbf{x}$ is dually residuated. Its dual residual is
equal to $$\Lambda_{\mathbf{a}}^\flat : \Is
\rightarrow \Is, \mathbf{x} \mapsto
\mathbf{a} \overline{\lslashh} \mathbf{x}= [\underline{a}\lslashh
\underline{x}, \underline{a}\lslashh \underline{x} \oplus
\overline{a} \lslashh \overline{x}].$$ Therefore, $ \mathbf{a}
\overline{\lslashh} \mathbf{b}$ is the smallest solution of
$\mathbf{a} \stackrel{-}{\odot} \mathbf{x} \succeq \mathbf{b} $,
and  equality is achieved if $\mathbf{b} \in \im \Lambda_{\mathbf{a}}$.
\end{proposition}
\begin{pf}
Let $ \Psi : \cod  \rightarrow \Is,
(\widetilde{x}',\widetilde{x}'') \mapsto
[\underline{x},\overline{x}]=[\widetilde{x}',\widetilde{x}'']$ be
the mapping which maps an ordered pair to an interval. This
mapping defines an isomorphism, since it is sufficient to deal with
the bounds to handle an interval. Then the result follows directly
from Proposition \ref{prop:dualresiduation_couple_ordonnes}.
 \qed
\end{pf}
\begin{cor}
Let $\sd$ be a semiring  and $\mathbf{A} \in \Is^{n\times p}$, $\mathbf{X}\in \Is^{p\times
q}$ and $\mathbf{Y}\in \Is^{n\times
q}$ be matrices. According to Corollary \ref{cor:dualresiduationmatrixproduct}, mapping $\Lambda_{\mathbf{A}} : \Is^{p \times q} \rightarrow
\Is^{n \times q},\mathbf{X} \mapsto \mathbf{A}
\stackrel{-}{\odot} \mathbf{X}$ is dually residuated. Its dual residual  is equal to
\begin{equation}\label{eq:dualresidutaionmatricesofinterval}
\Lambda_{\mathbf{A}}^\flat : \Is^{n \times q} \rightarrow
\Is^{p \times q}, \mathbf{Y} \mapsto
\mathbf{A}\overline{\lslashh}\mathbf{Y}= [\underline{A}\lslashh \underline{Y}, \underline{A}\lslashh \underline{Y} \oplus \overline{A}\lslashh \overline{Y}].
\end{equation}
\end{cor}
Additive closure and residuation being well defined over a semiring of intervals the Properties \ref{prop:etoilekleene} can be translated as follows.
\begin{propriete}\label{prop:intervaletoilekleene}
Let $\mathbf{A} \in \Is^{n \times n}$, $\mathbf{B} \in \Is^{n \times n}$,  $\mathbf{C} \in \Is^{n \times n}$  , and $\mathbf{X} \in \Is^{n \times p}$ be four matrices.  The following statements hold:
\begin{equation} \label{eq:intervalastarastarastar}
 \mathbf{A}^\ast \overline{\otimes} \mathbf{A}^\ast \overline{\otimes} \mathbf{X}=\mathbf{A}^\ast \overline{\otimes} \mathbf{X},
\end{equation}
\begin{equation} \label{eq:intervalastarsurastar}
\mathbf{A}^\ast \overline{\lslash} \mathbf{A}^\ast \overline{\lslash} \mathbf{X} = \mathbf{A}^\ast \overline{\lslash} \mathbf{X},
\end{equation}
\begin{equation} \label{eq:intervalastartimesastarsurastar}
\mathbf{A}^\ast \overline{\otimes} (\mathbf{A}^\ast \overline{\lslash} \mathbf{X})=\mathbf{A}^\ast \overline{\lslash} \mathbf{X},
\end{equation}
\begin{equation} \label{eq:intervalastarsurastartimesastar}
\mathbf{A}^\ast \overline{\lslash} (\mathbf{A}^\ast \overline{\otimes} \mathbf{X})=\mathbf{A}^\ast \overline{\otimes} \mathbf{X},
\end{equation}
\begin{equation} \label{eq:AstarGreaterBstarInterval}
 \mathbf{C}^\ast \preceq  \mathbf{A}^\ast \Leftrightarrow  \mathbf{A}^\ast   \mathbf{C}^\ast   \mathbf{X} =  \mathbf{A}^\ast  \mathbf{X} = \mathbf{C}^\ast \lslash ( \mathbf{A}^\ast  \mathbf{X}) \Leftrightarrow \im L_{ \mathbf{A}^\ast} \subset \im L_{ \mathbf{C}^\ast} \Leftrightarrow \im L_{ \mathbf{A}^\ast} \subset \im L_{ \mathbf{C}^\ast}^\sharp.
 \end{equation}
For the dual product the following property can be stated :
\begin{equation} \label{eq:intervalastarastarastar}
 \mathbf{B}_\ast \overline{\odot} \mathbf{B}_\ast \overline{\odot} \mathbf{X}=\mathbf{B}_\ast \overline{\odot} \mathbf{X},
\end{equation}
 and the following equivalences hold
\begin{equation} \label{eq:AstarXimAstarInterval}
\begin{array}{lclclclcl}
\mathbf{A} \overline{\otimes} \mathbf{X} \preceq \mathbf{X} & \Leftrightarrow & \mathbf{X}& =&\mathbf{A}^\ast \overline{\otimes} \mathbf{X}& \Leftrightarrow & \mathbf{A}^\ast \overline{\lslash} \mathbf{X} & \Leftrightarrow & \mathbf{X} \in \im L_{\mathbf{A}^\ast}, \\
\mathbf{X} \preceq \mathbf{B} \overline{\odot} \mathbf{X} & \Leftrightarrow &
\mathbf{X}&= &\mathbf{B}_\ast \overline{\odot} \mathbf{X} & \Leftrightarrow & \mathbf{B}_\ast \overline{\lslashh} \mathbf{X} & \Leftrightarrow & \mathbf{X} \in \im \Lambda_{\mathbf{B}_\ast}.
\end{array}
\end{equation}
\end{propriete}
\begin{rem}
According to Proposition \ref{res_interval} and \ref{dualres_interval}, the following implications hold :
\begin{equation}\nonumber
\begin{array}{lclclclcl}
\mathbf{X} \in \im L_{\mathbf{A}^\ast} \Rightarrow \mathbf{X}&=& [\underline{A}^\ast \underline{X},\overline{A}^\ast\overline{X}]  = [\underline{A}^\ast \lslash \underline{X} \wedge \overline{A}^\ast \lslash \overline{X},\overline{A}^\ast \lslash \overline{X}]\\
& =&[\underline{A}^\ast\lslash \underline{X},\overline{A}^\ast \lslash \overline{X}] \textnormal{ since }
\underline{A}^\ast\underline{X} \preceq\overline{A}^\ast \overline{X}, \\

\mathbf{X} \in \im \Lambda_{\mathbf{B}_\ast} \Rightarrow \mathbf{X}&=& [\underline{B}_\ast \odot \underline{X},\overline{B}_\ast\odot \overline{X}]  = [\underline{B}_\ast \lslashh \underline{X}, \underline{B}_\ast \lslash \underline{X} \oplus \overline{B}_\ast \lslashh \overline{X}]\\
& =&[\underline{B}_\ast\lslashh \underline{X},\overline{B}_\ast \lslashh \overline{X}] \textnormal{ since }
\underline{B}_\ast\odot \underline{X} \preceq\overline{B}_\ast \odot \overline{X}. \\
\end{array}
\end{equation}
\end{rem}

Below, the extension of Proposition \ref{propo:ProjectorInterBstarAstar}  to a  semiring of intervals is given.
\begin{proposition}\label{propo:IntervalProjectorInterBstarAstar}

Let $\sd$ be a semiring and $\mathbf{A}, \mathbf{B}  \in \Is^{n \times n}$ and $\mathbf{X} \in \sd^{n \times m}$.\\ If $\forall \mathbf{X}$  the equality $\mathbf{B_\ast} \overline{\lslashh} (\mathbf{A^\ast} \overline{\otimes} \mathbf{X})=(\mathbf{B_\ast} \overline{\lslashh} \mathbf{A^\ast}) \overline{\otimes} \mathbf{X}$ holds,  mapping
\begin{equation}\nonumber
\begin{array}{lclclcl}
\mathbf{P}  &:& \Is^{n \times m} & \rightarrow & \Is^{n \times m},
              \mathbf{X} & \mapsto &   (\mathbf{B_\ast}\overline{\lslashh} \mathbf{A}^\ast)^\ast \overline{\lslash} \mathbf{X}

\end{array}
\end{equation}
with \small
\begin{equation}\nonumber
\begin{array}{lclclcl}
 (\mathbf{B_\ast}\overline{\lslashh} \mathbf{A}^\ast)^\ast \overline{\lslash} \mathbf{X}& = &[ (( \underline{B}_\ast \lslashh \underline{A}^\ast)^\ast \lslash \underline{X}) \wedge
            (((\underline{B}_\ast \lslashh \underline{A}^\ast)\oplus( \overline{B}_\ast \lslashh \overline{A}^\ast))^\ast \lslash \overline{X}),
            (( \underline{B}_\ast \lslashh \underline{A}^\ast) \oplus ( \overline{B}_\ast \lslashh \overline{A}^\ast))^\ast\lslash \overline{X}
            ],
\end{array}
\end{equation}
\normalsize
is a projector in $\im L_{\mathbf{A^\ast}} \cap \im \Lambda_{\mathbf{B_\ast}}$, formally  $$\mathbf{P}(\mathbf{X})=\{ \bigvee \mathbf{Y} |\mathbf{Y}\preceq_{\Is} \mathbf{X} \text{ and } \mathbf{Y} \in \im L_{\mathbf{A^\ast}} \cap \im \Lambda_{\mathbf{B_\ast}} \}.$$

\end{proposition}
\begin{pf}
It is a direct application of Proposition \ref{propo:ProjectorInterBstarAstar}.
For the practical computation, from Proposition \ref{res_interval}, we get :
$$(\mathbf{B_\ast}\overline{\lslashh} \mathbf{A}^\ast)^\ast \overline{\lslash} \mathbf{X}=[
(( \underline{\mathbf{B}_\ast \overline{\lslashh} \mathbf{A}^\ast})^\ast \lslash \underline{X}) \wedge
            ( \overline{(\mathbf{B}_\ast \overline{\lslashh} \mathbf{A}^\ast})^\ast \lslash \overline{X}),
            ( \overline{(\mathbf{B}_\ast \overline{\lslashh} \mathbf{A}^\ast})^\ast \lslash \overline{X})
] $$
with, according to Propositions  \ref{dualres_interval} and \ref{prop:intervalclosuredualclosure},
$$
(\underline{\mathbf{B}_\ast \overline{\lslashh} \mathbf{A}^\ast})^\ast=(\underline{B}_\ast \lslashh \underline{A}^\ast)^\ast
$$
and
$$
(\overline{\mathbf{B}_\ast \overline{\lslashh} \mathbf{A}^\ast})^\ast =((\underline{B}_\ast \lslashh \underline{A}^\ast) \oplus (\overline{B}_\ast \lslashh \overline{A}^\ast))^\ast.
$$
\qed

\end{pf}

\begin{exmp}
Below, we compute the greatest interval vector which satisfies :
\begin{equation}
\nonumber
\begin{array}{lcl}
\mathbf{A}\overline{\otimes} \mathbf{X}\preceq \mathbf{X} \preceq \mathbf{B} \overline{\odot} \mathbf{X} \\
\mathbf{X}\preceq \mathbf{X_0},
\end{array}
\end{equation}
where
\small
$$\mathbf{A}=\begin{pmatrix} [\varepsilon, \varepsilon] & [\varepsilon,\varepsilon] &[\varepsilon,\varepsilon] & [\varepsilon,\varepsilon] & [\varepsilon,\varepsilon] \\
[7,11] & [\varepsilon,\varepsilon] &[8,14] & [\varepsilon,\varepsilon] & [2,7] \\
[\varepsilon,\varepsilon] & [\varepsilon,\varepsilon] &[\varepsilon,\varepsilon] & [\varepsilon,\varepsilon] & [\varepsilon,\varepsilon] \\
 [\varepsilon,\varepsilon] & [\varepsilon,\varepsilon] &[4,12] & [\varepsilon,\varepsilon] & [1,5] \\
[\varepsilon,\varepsilon] & [\varepsilon,\varepsilon] &[\varepsilon,\varepsilon] & [\varepsilon,\varepsilon] & [\varepsilon,\varepsilon]
 \end{pmatrix}, \text{ }\mathbf{B}= \begin{pmatrix} [\top,\top] & [\top,\top] &[\top,\top] & [\top,\top] & [\top,\top] \\
[11,16] & [\top,\top]&  [15,19] & [\top,\top] & [7,10] \\
[\top,\top] & [\top,\top] &[\top,\top] & [\top,\top] & [\top,\top] \\
[\top,\top] & [\top,\top] &[13,18] & [\top,\top] & [5,9] \\
[\top,\top] & [\top,\top] &[\top,\top] & [\top,\top] & [\top,\top] \end{pmatrix}$$
and $\mathbf{X_0}=\begin{pmatrix} [10,14] & [10,14] & [10,14] & [10,14]& [10,14] \end{pmatrix}^T.$ \\
\normalsize
We get :

 $$(\underline{B}_\ast \lslashh \underline{A}^\ast)^\ast=\begin{pmatrix} e & -11 & -3 & -14 & -9 \\
7 & e & 8 & -3 & 2 \\
 -8& -15 & e & -13 & -12 \\
1 & -6 & 4 & e & 1 \\
e & -7 & 1 & -5 & e
 \end{pmatrix}, \text{ }(\overline{B}_\ast \lslashh \overline{A}^\ast)^\ast=\begin{pmatrix} e & -16 & -2 & -18 & -9 \\
11 & e & 14 & -2 & 7 \\
 -8& -19 & e & -18 & -12 \\
6 & -5 & 12 & e & 5 \\
1 & -10 & 4 & -9 & e
 \end{pmatrix}.$$
This yields $\mathbf{X}=\mathbf{P(X_0)}=\begin{pmatrix}
 [3,3] & [10,14] & [0,0] & [10,12] & [7,7] \end{pmatrix}^T$ as greatest interval vector.
 \end{exmp}
 \begin{exmp}
We provide also an example in the semiring $\zmaxgamma$. We consider :
\small
$$\mathbf{A}=\begin{pmatrix} [\varepsilon, \varepsilon] & [\varepsilon,\varepsilon] &[8\gamma^2,8\gamma] \\
[\varepsilon, \varepsilon] & [\varepsilon,\varepsilon] &[\varepsilon,\varepsilon] \\
[7\gamma \oplus 9 \gamma^2,10 \oplus 11 \gamma^3] & [2\gamma \oplus 4 \gamma^3, 4 \gamma \oplus 6 \gamma^2] & [\varepsilon,\varepsilon] \\
 \end{pmatrix}, \text{ }\mathbf{B}= \begin{pmatrix} [\top,\top] & [\top,\top] &[15\gamma,18] \\
 [\top,\top] & [\top,\top] &[\top,\top]\\
 [\top,\top] & [5\gamma,7] &[\top,\top]
 \end{pmatrix}
 $$

 and $\mathbf{X_0}=\begin{pmatrix} [4\gamma \oplus 7 \gamma^4 (18 \gamma)^\ast,7 \oplus 8 \gamma^3 (18 \gamma)^\ast]\\ [5\gamma^2 \oplus 8 \gamma^5 (18 \gamma)^\ast,8 \gamma \oplus 9 \gamma^4 (18 \gamma)^\ast]\\ [6 \gamma^3 \oplus9 \gamma^6 (18 \gamma)^\ast,9\gamma^2 \oplus 10 \gamma^5 (18 \gamma)^\ast] \end{pmatrix}.$ \\
\normalsize
According to the computation rules given in Definition \ref{def:zmaxgamma} (see also \cite{gaubert92a,lhommeau00a} for algorithmic issues and software tools), the following vector is obtained :
 $$\mathbf{X}=\mathbf{P(X_0)}=\begin{pmatrix} [21\gamma^4(18\gamma)^\ast,17\gamma^3(18\gamma)^\ast]\\
 [4\gamma^2(18\gamma)^\ast,5\gamma(18\gamma)^\ast]\\
 [6\gamma^3(18\gamma)^\ast,9\gamma^2(18\gamma)^\ast]
 \end{pmatrix}.$$

 \end{exmp}

\section{Conclusion}
This work deals with a dual product in a semiring and its extension to  semirings of intervals.
Sufficient conditions are given in order to ensure the existence of a
projector in the solution set of the following system :
$\mathbf{A} \overline{\otimes} \mathbf{X} \preceq \mathbf{X} \preceq \mathbf{B} \overline{\odot} \mathbf{X}$, where $\mathbf{A}$, $\mathbf{B}$ and $\mathbf{X}$ are interval matrices.
This projector can be  useful to solve control problems for timed discrete event systems.
More precisely, control for uncertain systems with parameters that are only known to be in an interval,
and where the state evolution is subject to constraints (see e.g. \cite{ouerghi06a,ouerghi06b,katz07,hardouin10a,maia11a,brunsch11}).

\bibliographystyle{plain}        

\end{document}